\documentclass[12pt]{amsart}
\newtheorem{thm}{Theorem}[section]
\newtheorem{theorem}[thm]{Theorem}
\newtheorem{corollary}[thm]{Corollary}
\newtheorem{lemma}[thm]{Lemma}
\newtheorem{proposition}[thm]{Proposition}
\theoremstyle{definition}
\newtheorem{definition}[thm]{Definition}
\newtheorem{example}[thm]{Example}

\theoremstyle{remark}
\newtheorem{remark}[thm]{Remark}
\numberwithin{equation}{section}

\newcommand{\brm}{\mathrm{b}}
\newcommand{\crm}{\mathrm{c}}
\newcommand{\erm}{\mathrm{e}}
\newcommand{\Irm}{\mathrm{I}}
\newcommand{\mrm}{\mathrm{m}}
\newcommand{\prm}{\mathrm{p}}

\newcommand{\wrm}{\mathrm{w}}

\newcommand{\DC}{\mathcal{D}}
\newcommand{\GC}{\mathcal{G}}
\newcommand{\PC}{\mathcal{P}}
\newcommand{\SC}{\mathcal{S}}
\newcommand{\UC}{\mathcal{U}}
\newcommand{\VC}{\mathcal{V}}
\newcommand{\WC}{\mathcal{W}}

\newcommand{\Gb}{\mathbf{G}}
\newcommand{\Jb}{\mathbf{J}}
\newcommand{\ub}{\mathbf{u}}
\newcommand{\wb}{\mathbf{w}}
\newcommand{\xb}{\mathbf{x}}
\newcommand{\yb}{\mathbf{y}}
\newcommand{\zb}{\mathbf{z}}

\newcommand{\R}{\mathbb R}
\newcommand{\C}{\mathbb C}

\newcommand{\Sf}{{\mathfrak S}}
\newcommand{\Tf}{{\mathfrak T}}

\newcommand{\al}{\alpha}
\newcommand{\de}{\delta}
\newcommand{\ve}{\varepsilon}
\newcommand{\F}{\Phi}

\newcommand{\la}{\lambda}
\newcommand{\La}{\Lambda}
\newcommand{\Om}{\Omega}

\newcommand{\ps}{\psi}
\newcommand{\si}{\sigma}
\newcommand{\Si}{\Sigma}
\renewcommand{\th}{\theta}
\newcommand{\The}{\Theta}

\newcommand{\rank}{\mathrm{rank}\,}
\newcommand{\ess}{\mathrm{ess}}

\newcommand{\ch}{\mathrm{conv}\,}
\newcommand{\cch}{\overline{\mathrm{conv}}\,}
\newcommand{\ex}{\mathrm{ex}\,}
\newcommand{\sym}{\mathrm{sym}}

\begin{document}

\title[Birkhoff's theorem]
{Birkhoff's theorem and multidimensional numerical range}
\author{Yu. Safarov}

\address{Department of Mathematics, King's College,
Strand, London WC2R 2LS} \email{ysafarov@mth.kcl.ac.uk}



\date{August 2004}


\begin{abstract}
We show that, under certain conditions, Birkhoff's theorem on
doubly stochastic matrices remains valid for countable families of
discrete probability spaces which have nonempty intersections.
Using this result, we study the relation between the spectrum of a
self-adjoint operator $\,A\,$ and its multidimensional numerical
range. It turns out that the multidimensional numerical range is a
convex set whose extreme points are sequences of eigenvalues of
the operator $\,A\,$. Every collection of eigenvalues which can be
obtained by the Rayleigh--Ritz formula generates an extreme point
of the multidimensional numerical range. However, it may also have
other extreme points.
\end{abstract}

\maketitle

Recall that a (possibly infinite) matrix is said to be {\sl doubly
stochastic} if all its entries are non-negative and the sum of
entries in every row and every column is equal to one. Birkhoff's
theorem \cite{B} says that
\begin{enumerate}
\item[{\bf(i)}]
the extreme points of the convex set of doubly
stochastic matrices are permutation matrices and
\item[{\bf(ii)}]
the set of doubly stochastic matrices coincides with the closed
convex hull of the set of permutation matrices.
\end{enumerate}
The first aim of this paper is to show that, under certain
conditions, Birkhoff's theorem remains valid for a countable
family of discrete probability spaces which have nonempty
intersections (see Remark \ref{R2.1}). We join every two points
lying in the same probability space by an edge and reformulate the
problem in terms of weighted graphs. It turns out that {\bf(i)}
and {\bf(ii)} hold true whenever the underlying graph satisfies
the conditions ({\bf g}$_1$)--({\bf g}$_3$) introduced in Section
2. The conditions ({\bf g}$_1$) and ({\bf g}$_3$) are purely
technical and can probably be removed or weakened. The geometric
condition ({\bf g}$_2$) is necessary (see Remark \ref{R2.5}).

The second aim of the paper is to study the relation between the
spectrum of a self-adjoint operator $\,A\,$ and its
$\,m$-dimensional numerical range $\,\Si(m,A)\,$. The latter is
defined as the set of all $\,m$-dimensional vectors of the form
$\,\{Q_A[u_1],Q_A[u_2],\ldots\}\,$, where $\,Q_A\,$ is the
corresponding quadratic form,
$\,\{u_1,u_2,\ldots\}\subset\DC(Q_A)\,$ is an arbitrary
orthonormal set containing $\,m\,$ elements and
$\,m=1,2,\ldots,\infty\,$. Using an infinite dimensional version
of Birkhoff's theorem, we prove that
\begin{enumerate}
\item[{\bf(1)}]
the $\,m$-dimensional numerical range $\,\Si(m,A)\,$ is a convex
set,
\item[{\bf(2)}]
the extreme points of $\,\Si(m,A)\,$ belong to the corresponding
$\,m$-dimen\-sional point spectrum $\,\si_p(m,A)\,$,
\item[{\bf(3)}]
every collection of $\,m\,$ lowest or highest eigenvalues which
can be found with the use of the Rayleigh--Ritz formula generates
an extreme point of $\,\Si(m,A)\,$,
\item[{\bf(4)}]
the extreme points of the closure $\,\overline{\Si(m,A)}\,$ belong
to the $\,m$-dimen\-sional spectrum $\,\si(m,A)\,$,
\item[{\bf(5)}]
the closed convex hull of $\,\si(m,A)\,$ coincides with
$\,\overline{\Si(m,A)}\,$
\end{enumerate}
(see Section 4 for precise statements and definitions). The item
{\bf(3)} can be regarded as a geometric version of the variational
principle. The set $\,\Si(m,A)\,$ may also have other extreme
points (see Remark \ref{R4.12}). Therefore one can obtain more
information about the point spectrum by studying the extreme
points of $\,\Si(m,A)\,$ than by applying the standard variational
formulae.

The paper is organised as follows. For the sake of convenience, in
Section 1 we give definitions and results on sequence spaces and
locally convex topologies, which are used throughout other
sections. Almost all these results are well known; most of them
can be found in \cite{K}, Sections 20.9, 21.2, 30 and \cite{Ru},
Section 2.4.

Section 2 is devoted to Birkhoff's theorem. Many proofs of this
theorem are known for finite matrices (see, for example, \cite{MO}
or \cite{BP}). The problem of extending {\bf(i)} and {\bf(ii)} to
infinite matrices is known as Birkhoff's problem $111$. It has
been studied in \cite{Gr}, \cite{Is}, \cite{Ke}, \cite{Le},
\cite{Mu} and \cite{RP}. However, their results are not sufficient
for our purposes because
\begin{enumerate}
\item[(i)]
in order to deal with unbounded operators, we need {\bf(i)} not
only for the whole set of stochastic matrices but also for some
its subsets which were not considered in these papers,
\item[(ii)]
we need {\bf(ii)} with respect to a finer topology than the
topology introduced in \cite{Ke} or \cite{RP}, whereas \cite{Is}
deals with a too strong topology such that {\bf(ii)} does not hold
true.
\end{enumerate}

Our proof of {\bf(i)} and {\bf(ii)} is based on the well known
idea of shifting weights along edges of the underlying graph. It
is almost purely combinatorial and works equally well for finite
and infinite weighted graphs or matrices. Formally speaking, in
Sections 3 and 4 we consider only infinite matrices. However, in
the proof of Theorem \ref{T3.15} we apply results related to more
general weighted graphs. For infinite graphs and matrices
{\bf(ii)} depends upon the choice of an appropriate topology. We
give an explicit description of the strong and Mackey topologies
on the set of (sub)stochastic weights (Corollaries \ref{C2.11} and
\ref{C2.12}), and show that {\bf(ii)} holds true with respect to
the Mackey topology (Theorem \ref{T2.15}), but not necessarily
with respect to the strong topology (Example \ref{E2.19}).

In Section 3 we consider operators generated by stochastic
matrices and derive a number of corollaries from Birkhoff's
theorem. Many of these results seem almost obvious. However, our
proofs of the key Theorems \ref{T3.10} and \ref{T3.15} are
surprisingly long and complicated. It is not clear whether they
can be essentially simplified.

Section 4 is about multidimensional spectra and numerical ranges.
Here we give precise statements and proofs of {\bf(1)}--{\bf(5)}
for a self-adjoint operator $\,A\,$ (see Corollaries \ref{C4.7},
\ref{C4.11} and Lemma \ref{L4.10}). The corresponding results for
finite matrices $\,A\,$ are well known and rather elementary (see,
for example, \cite{AU} or \cite{MO}). If $\,A\,$ is compact, one
can probably obtain {\bf(1)}--{\bf(5)} by considering its finite
dimensional approximations (in \cite{Ma1} and \cite{Si} similar
ideas have been used for studying $s$-numbers of compact
operators). However, the general case is much more complex as the
operator $\,A\,$ may have continuous spectrum or (and) several
accumulation points of its discrete spectrum, which makes it
impossible to find an effective approximation procedure. In the
end of Section 4 we prove two variational formulae (Corollaries
\ref{C4.16} and \ref{C4.17}) and show that $\,\si(m,A)\,$ is a
subset of the closed convex hull of $\,\bigcup_\th\si(m,A_\th)\,$
whenever the self-adjoint operator $\,A\,$ belongs to the closed
convex hull of the family of self-adjoint operators $\,A_\th\,$
(Corollary \ref{C4.22}); all these results are simple consequences
of {\bf(1)}--{\bf(5)}.

There are many other concepts of multidimensional numerical range
\cite{BD}, \cite{H}, \cite{LMMT}. We briefly discuss some of them
in Subsection \ref{S4.1}.

\medskip \noindent {\bf Acknowledgements.} I would like to thank
E.B. Davies, A. Markus, Yu. Netrusov and M. Solomyak for their
valuable comments and constructive criticism. I am also grateful
to L. Landau for his useful remarks on Birkhoff's theorem.

\section{Sequence spaces}

\subsection{Notation and definitions}
Let
\begin{enumerate}
\item[]
$\,\hat\R:=[-\infty,+\infty]\,$,
\item[]
$\R^\infty$ be the linear spaces of all real sequences
$\,\xb=\{x_1,x_2,\ldots\}\,$,
\item[]
$\R_0^\infty$ be the subspace of sequences which converge to zero
and
\item[]
$\R_{00}^\infty$ be the subspace of sequences with finitely many
nonzero entries.
\end{enumerate}
We shall often consider the Euclidean space $\R^m$ as a finite
dimensional subspace of $\,\R_{00}^\infty\,$, so that the
$\,m$-dimensional real vector $\,(x_1,x_2,\ldots,x_m)\,$ is
identified with the sequence
$\,(x_1,x_2,\ldots,x_m,0,0,\ldots)\,$. If $\,\xb\in\R^\infty\,$,
let
\begin{align}\label{1.1}
|\xb|\ \ &:=\ \{|x_1|,|x_2|,\ldots\}\,,\nonumber\\
\xb^{(m)}\ &:=\
\{x_1,x_2,\ldots,x_m,0,0,\ldots\}\,,\qquad m=1,2,\ldots,\\
\xb^{(\infty)}\ &:=\ \xb\,.\nonumber
\end{align}

Throughout the paper $X$ denotes a real linear subspace of
$\,\R^\infty\,$ endowed with a locally convex topology $\,\Tf\,$
and $X^*$ is its dual space. We shall always be assuming that
$\,\Tf\,$ is finer (that is, not weaker) than the topology of
element-wise convergence.

If $\,\Om\,$ is a subset of $X$ then $\,\ex\Om\,$, $\,\ch\Om\,$,
$\,\cch\Om\,$ denote the set of extreme points of $\,\Om\,$, the
convex hull of $\,\Om\,$ and its $\,\Tf$-closure respectively.
Recall that $\,\xb\in\Om\,$ is called an {\sl extreme point} of
$\,\Om\,$ if $\,\xb\,$ cannot be represented as a convex linear
combination of two other elements of $\,\Om\,$. If the set
$\,\Om\,$ is $\,\Tf$-compact then, according to the Krein--Milman
theorem, $\,\cch\Om=\cch(\ex\Om)\,$. An element $\,\xb\in\Om\,$ is
said to be {\sl $\,\Tf$-exposed} if there exists a linear
$\,\Tf$-continuous functional $\,\xb^*\in X^*\,$ such that
$\,\langle\xb,\xb^*\rangle>\langle\yb,\xb^*\rangle\,$ for all
$\,\yb\in\Om\,$. Every exposed point of $\,\Om\,$ belongs to
$\,\ex\Om\,$ but an extreme point is not necessarily exposed.

Denote by $\,X'\,$ the linear space of all real sequences
$\,\xb'=\{x'_1,x'_2,\ldots\}\in\R^\infty\,$ such that
$\,\sum_{i=1}^\infty|x_i\,x'_i|<\infty\,$ for all $\,\xb\in X\,$.
If $\,X''=X\,$ then the space $X$ is said to be {\sl perfect\/}.
We have $\,X\subseteq X''\,$ and $\,\R_{00}^\infty\subseteq
X'=X'''\,$; in particular, $X'$ is perfect. The intersection of an
arbitrary collection of perfect spaces is perfect. However, the
linear span of a collection of perfect spaces may not be perfect.
For example, if $\,X\,$ is a one dimensional subspace of
$\,\R_0^\infty\,$ then $\,X''\subset\R_0^\infty\,$ but
$\,(\R_0^\infty)''=l^\infty\,$.

The set of sequences $\,\tilde\xb=\{\tilde x_1,\tilde
x_2,\ldots\}\,$ such that $\,|\tilde x_j|\le|x_j|\,$ for some
$\,\xb\in\Om\,$ is said to be the {\sl normal cover\/} of the set
$\,\Om\,$. A set (or subspace) of $\,\R^\infty\,$ is said to be
{\sl normal\/} if it coincides with its normal cover. We have
$\,X'=(\tilde X)'\,$, where $\,\tilde X\,$ is a normal cover of
$\,X\,$. Therefore a perfect space is normal.

\subsection{Topologies on sequence spaces}
Every sequence $\,\xb'\in X'\,$ defines the linear functional
$\,\langle\xb,\xb'\rangle:=\sum_{j=1}^\infty x_j\,x'_j\,$ on the
space $\,X\,$. Further on we shall always be assuming that
$\,\R_{00}^\infty\subseteq X\,$. Then every nonzero element of
$\,X'\,$ defines a nonzero functional and therefore we can
introduce the weak${}^*$ topology $\,\Tf_\wrm(X',X)\,$ on
$\,X'\,$. If $\,\Sf\,$ is an arbitrary family of weak${}^*$
bounded sets $\,\Om'\in X'\,$ then the family of seminorms
\begin{equation}\label{1.2}
p_{\Om'}(\xb)\ :=\
\sup_{\xb'\in\Om'}|\langle\xb,\xb'\rangle|\,,\qquad\Om'\in\Sf\,,
\end{equation}
defines a locally convex topology on the space $\,X\,$, which is
usually called the $\,\Sf$-topology. We shall deal with the
following $\,\Sf$-topologies on $\,X\,$:
\begin{enumerate}
\item[(1)]
the topology of element-wise convergence $\,\Tf_0\,$, generated by
the family $\,\Sf\,$ of all finite subsets of
$\,\R_{00}^\infty\,$;
\item[(2)]
the weak topology $\,\Tf_\wrm(X,X')\,$, generated by the family
$\,\Sf\,$ of all finite subsets of $\,X'\,$;
\item[(3)]
the Mackey topology $\,\Tf_\mrm(X,X')\,$, generated by the family
$\,\Sf\,$ of all absolutely convex $\,\Tf_\wrm(X',X)$-compact
subsets of $\,X'\,$;
\item[(4)]
the strong topology $\,\Tf_\brm(X,X')\,$, generated by the family
$\,\Sf\,$ of all $\,\Tf_\wrm(X',X)$-bounded subsets of $\,X'\,$.
\end{enumerate}

Every next topology in this list is finer than the previous one.
Each of them is equivalent to the usual Euclidean topology
whenever $\,\dim X<\infty\,$. The strong topology
$\,\Tf_\brm(X,X')\,$ is generated by all lower
$\,\Tf_\wrm(X,X')$-semi\-continuous seminorms on $\,X\,$ and the
Mackey topology $\,\Tf_\mrm(X,X')\,$ is defined by all lower
$\,\Tf_\wrm(X,X')$-semi\-continuous seminorms $\,p\,$ on $\,X\,$
such that
\begin{equation}\label{1.3}
p(\xb-\xb^{(m)})\ \underset{m\to\infty}\to\ 0\,,\qquad
\forall\xb\in X\,.
\end{equation}

The perfect space $\,X''\,$ is obtained from $\,X\,$ by adding all
$\,\Tf_0$-limits of $\,\Tf_\wrm(X,X')$-Cauchy sequences in
$\,X\,$. A perfect space $\,X\,$ is $\,\Tf_\brm(X,X')$-comp\-lete,
$\,\Tf_\mrm(X,X')$-complete and sequentially
$\,\Tf_\wrm(X,X')$-complete but is not necessarily
$\,\Tf_\wrm(X,X')$-complete. By the Mackey--Arens theorem,
$\,\Tf_\mrm(X,X')\,$ is the finest locally convex topology on the
space $\,X\,$ such that its topological dual $\,X^*\,$ coincides
with $\,X'\,$. If $\,X'\,$ is $\,\Tf_\wrm(X',X)$-complete then the
$\,\Tf_\brm(X,X')$-dual of $\,X\,$ also coincides with $\,X'\,$.

By Mackey's theorem, a subset of a locally convex space is weakly
bounded if and only if it is bounded in any topology
generating the same dual space. For a sequence space $\,X\,$, we
have the following stronger result
which implies that $\,\Om\subset X\,$ is
$\,\Tf_\wrm(X,X')$-bounded if and only if it is
$\,\Tf_\brm(X,X')$-bounded.

\begin{theorem}\label{T1.1}
Assume that $\,\Om\subset X\,$ is $\,\Tf_\wrm(X,X')$-bounded and
$\,\Om'\subset X'\,$ is $\,\Tf_\wrm(X',X)$-bounded. Then the set
of sequences $\,\{x_1\,x'_1\,,x_2\,x'_2\,,\ldots\}\,$, where
$\,\xb=\{x_1,x_2,\ldots\}\in\Om\,$ and
$\,\xb'=\{x'_1,x'_2,\ldots\}\in\Om'\,$, is bounded in $\,l^1\,$.
\end{theorem}

\begin{proof}
See \cite{Ru}, Chapter 2, Proposition 1.4.
\end{proof}

The following theorem can be proved in the same way as Theorem 2.4
in \cite{Ru}, Chapter 2, where the author assumed that $\,X\,$ is
perfect.

\begin{theorem}\label{T1.2}
If $\,X\,$ is a normal space and $\,\Om'\subset X'\,$ then the
following two conditions are equivalent:
\begin{enumerate}
\item[(1)]
$\,\Om'\,$ is $\,\Tf_\wrm(X',X)$-compact,
\item[(2)]
$\,\Om'\,$ is $\,\Tf_0$-compact and
$\,\lim\limits_{n\to\infty}\,\sup\limits_{\xb'\in\Om'}\,
\sum\limits_{i=n}^\infty|x_i\,x'_i|=0\,$ for each $\,\xb\in X\,$.
\end{enumerate}
\end{theorem}

\begin{remark}\label{R1.3}
If $\,\{\xb_n\}\subset X\,$ is a $\,\Tf_\wrm(X,X')$-Cauchy
sequence which converges to $\,\xb\in X''\,$ in the topology
$\,\Tf_0\,$, then by Fatou's lemma
$$
\sup_{\xb'\in\Om'}|\langle\xb,\xb'\rangle|\ \le\
\sup_{\xb'\in\Om'}\langle|\xb|,|\xb'|\rangle\ \le\
\sup_{\xb'\in\Om'}\sup_n\langle|\xb_n|,|\xb'|\rangle\,.
$$
Since the Cauchy sequence $\,\{\xb_n\}\,$ is
$\,\Tf_\wrm(X,X')$-bounded, Theorem \ref{T1.1} and the above
inequality imply that the set $\,\Om'\subset X'\,$ is
$\,\Tf_\wrm(X',X'')$-bounded if and only if it is
$\,\Tf_\wrm(X',X)$-bounded. Therefore the strong topology
$\,\Tf_\brm(X,X')\,$ coincides with the restriction of
$\,\Tf_\brm(X'',X')\,$ to $\,X\,$. However, this is not
necessarily the case with the Mackey topologies.
\end{remark}

\begin{example}\label{E1.4}
If $\,X=\R_0^\infty\,$ then $\,X'=l^1\,$, $\,X''=l^\infty\,$ and
$\,\Tf_\brm(l^\infty,l^1)\,$ is the $\,l^\infty$-topology. Theorem
\ref{T1.2} implies that the closed unit ball in the space
$\,l^1\,$ is $\,\Tf_\wrm(l^1,\R_0^\infty)$-compact. Therefore
$\,\Tf_\mrm(\R_0^\infty,l^1)=
\left.\Tf_\brm(l^\infty,l^1)\right|_{\R_0^\infty}\,$. The Mackey
topology $\,\Tf_\mrm(l^\infty,l^1)\,$ on $\,\R_0^\infty\,$ is
strictly coarser than $\,\Tf_\mrm(\R_0^\infty,l^1)\,$. Indeed, if
$\,\xb=\{1,1,\ldots\}\,$ and
$\,\tilde\xb_m:=\xb^{(m+1)}-\xb^{(m)}\,$ then
$\,\tilde\xb_m\in\R_0^\infty\,$,
$\,\|\tilde\xb_m\|_{l^\infty}=1\,$ but, by Theorem \ref{T1.2},
$\,\tilde\xb_m\to0\,$ as $\,m\to\infty\,$ in the topology
$\,\Tf_\mrm(l^\infty,l^1)\,$.
\end{example}

\begin{remark}\label{R1.5}
Let $\,\tilde\Om'\,$ be the normal cover of the set $\,\Om'\subset
X'\,$. Theorem \ref{T1.1} implies that $\,\tilde\Om'\,$ is
$\,\Tf_\wrm(X',X)$-bounded whenever $\,\Om'\,$ is
$\,\Tf_\wrm(X',X)$-bounded. If $\,X\,$ is normal then, by Theorem
\ref{T1.2}, $\,\tilde\Om'\,$ is $\,\Tf_\wrm(X',X)$-compact
whenever $\,\Om'\,$ is $\,\Tf_\wrm(X',X)$-compact. Obviously,
$$
p_{\Om'}(\xb)=\sup_{\xb'\in\Om'}|\langle\xb,\xb'\rangle|
\le\sup_{\xb'\in\tilde\Om'}|\langle\xb,\xb'\rangle|
=\sup_{\xb'\in\tilde\Om'}\sum_{j=1}^\infty|x_j|\,|x'_j|
=p_{\tilde\Om'}(\xb)
$$
and the seminorms $\,p_{\tilde\Om'}\,$ are lower
$\,\Tf_0$-semi\-continuous. Therefore the strong topology
$\,\Tf_\brm(X,X')\,$ on an arbitrary space $\,X\,$ is generated by
all lower $\,\Tf_0$-semi\-continuous seminorms and the Mackey
topology $\,\Tf_\mrm(X,X')\,$ on a normal space $\,X\,$ is
generated by all lower $\,\Tf_0$-semi\-continuous seminorms
satisfying (\ref{1.3}).
\end{remark}

\subsection{Symmetric sequence spaces}
Our choice of notation in the following definition will become
clear in Section 3.

\begin{definition}\label{D1.6}
If $\,\xb\in\R^\infty\,$, let
\begin{enumerate}
\item[]
$P_\xb\,$ be the set of all sequences $\yb\in\R^\infty$ obtained
from the sequence $\,\xb\,$ by permutations of its entries,
\item[]
$P_\xb^r\,$ be the set of all sequences $\tilde\yb\in\R^\infty$
whose entries form a subsequence of a sequence $\yb\in P_\xb$ and
\item[]
$\,P_\xb^\emptyset\,$ be the set of all sequences obtained from
sequences $\tilde\yb\in P_\xb^r\,$ by adding an arbitrary
collection of zero entries.
\end{enumerate}
We shall say that a sequence space $\,X\,$ is {\sl symmetric\/} if
$\,P_\xb\subset X\,$ for every $\,\xb\in X\,$. A seminorm $\,p\,$
on a symmetric space $\,X\,$ is said to be symmetric if
$\,p(\yb)=p(\xb)\,$ whenever $\,\yb\in P_\xb\,$.
\end{definition}

If $\,X\,$ is symmetric then $\,P_{\xb'}^\emptyset\subset X'\,$
for every $\,\xb'\in X'\,$. The seminorm $\,p_{\Om'}\,$ defined by
(\ref{1.2}) is symmetric if and only if
$\,\Om'=\bigcup_{\xb'\in\Om'}P_{\xb'}\,$. The following result is
a consequence of Theorems \ref{T1.1} and \ref{T1.2} (see Remark
\ref{R3.2}).

\begin{corollary}\label{C1.7}
Let $\,X\,$ be a symmetric space such that
$\,X\not\subset\R_{00}^\infty\,$, $\,\Om'\,$ be a subset of
$\,X'\,$ and
$\,\Om'_\sym:=\bigcup_{\xb'\in\Om'}P_{\xb'}^\emptyset\,$. If
$\,X\subseteq l^\infty\,$ then $\,\Om'_\sym\,$ is
$\,\Tf_\wrm(X',X)$-bounded whenever $\,\Om'\,$ is
$\,\Tf_\wrm(X',X)$-bounded. If $\,X\subseteq\R_0^\infty\,$ and
$\,X\,$ is normal then $\,\Om'_\sym\,$ is
$\,\Tf_\wrm(X',X)$-compact whenever $\,\Om'\,$ is
$\,\Tf_\wrm(X',X)$-compact.
\end{corollary}

By Corollary \ref{C1.7}, if $\,X\,$ is a symmetric subspace of
$\,l^\infty\,$ and $\,X\not\subset\R_{00}^\infty\,$ then the
strong topology $\,\Tf_\brm(X,X')\,$ is generated by a family of
symmetric $\,\Tf_0$-semi\-continuous seminorms $\,p\,$ such that
\begin{equation}\label{1.4}
p(\yb)\ \le\ p(\xb)\,,\qquad\forall\yb\in X\bigcap
P_\xb^\emptyset\,,\quad\forall\xb\in X\,.
\end{equation}
If $\,X\,$ is a normal symmetric subspace of $\,\R_0^\infty\,$ and
$\,X\not\subset\R_{00}^\infty\,$ then the Mackey topology
$\,\Tf_\mrm(X,X')\,$ is generated by a family of symmetric
$\,\Tf_0$-semi\-continuous seminorms $\,p\,$ satisfying
(\ref{1.3}) and (\ref{1.4}).

\begin{example}\label{E1.8}
If $\,X=\R^\infty\,$ then $\,X'=\R_{00}^\infty\,$ and
$\,\Tf_0=\Tf_\mrm(\R^\infty,\R_{00}^\infty)
=\Tf_\brm(\R^\infty,\R_{00}^\infty)\,$.
This topology cannot be defined with the use of symmetric
seminorms. If $\,X=l^p\,$ with $1\le p\le\infty\,$ then
$\,X'=l^{p'}\,$ and $\,\Tf_\brm(l^p,l^{p'})\,$ is the usual
$\,l^p$-topology. If $\,p<\infty\,$ then
$\,\Tf_\brm(l^p,l^{p'})=\Tf_\mrm(l^p,l^{p'})\,$, but the Mackey
topology $\,\Tf_\mrm(l^\infty,l^1)\,$ is strictly coarser than the
$\,l^\infty$-topology and is not generated by a family of
symmetric seminorms.
\end{example}

\begin{example}\label{E1.9}
Let $\,\F\,$ be a symmetric lower $\,\Tf_0$-semicontinuous
Schatten norm on $\,\R_0^\infty\,$ and $\,s_\F^{(0)}\subseteq
s_\F\subseteq l^\infty\,$ be the corresponding linear subspaces of
sequences (see, for example, \cite{Si} or \cite{Ma1}; in the
latter paper $\,\F\,$ is called a symmetric gauge function and the
corresponding subspaces are denoted by $\,l_\F\,$ and
$\,l_\F^{(0)}\,$). Then the norm topology on a subspace
$\,X\subset s_\F\,$ is always coarser than $\,\Tf_\brm(X,X')\,$
and is coarser than $\,\Tf_\mrm(X,X')\,$ whenever $\,X\subset
s_\F^{(0)}\,$.
\end{example}

\begin{example}\label{E1.10}
Let $\,\xb\in\R_0^\infty\,$, $\,\xb\not\in l^1\,$ and $\,X\,$ be
the subspace spanned by the normal cover $\,\tilde P_\xb\,$ of the
set $\,P_\xb\,$. Then $\,X'\,$ consists of all sequences
$\,\xb'\in\R_0^\infty\,$ such that
\begin{equation}\label{1.5}
\|\xb'\|_L\ :=\ \sup_{\yb\in\tilde
P_\xb}|\langle\yb,\xb'\rangle|\ <\ \infty\,.
\end{equation}
The space $\,X'\,$ provided with the norm (\ref{1.5}) is called
the {\sl Lorentz} space associated with the weight sequence
$\,\xb\,$ (see, for example, \cite{LT}, Section 4.e). We have
\begin{equation}\label{1.6}
\sum_{k=1}^\infty |y_j|\,|x'_j|^*\ =\
\sum_{m=1}^\infty(|x'_m|^*-|x'_{m+1}|^*)\,\sum_{j=1}^m|y_j|\,,
\qquad\forall\xb',\yb\in\R_0^\infty\,,
\end{equation}
where $\,\{|z_1|^*,|z_2|^*\,\ldots\}\,$ denotes either the
non-increasing rearrangement of the sequence $\,|\zb|\,$ or (if
$\,|\zb|\,$ contains infinitely many nonzero entries and at least
one zero entry) the non-increasing rearrangement of its nonzero
entries. Using this identity, one can easily show that $\,\yb\in
X''\,$ if and only if
\begin{equation}\label{1.7}
\|\yb\|_M\ :=\ \sup_{m\ge1}\,R_m(|\yb|)\,(R_m(|\xb|))^{-1}\ =\
\sup_{\|\xb'\|_L<1}|\langle\yb,\xb'\rangle| \ <\ \infty\,,
\end{equation}
where $\,R_m(|\zb|):=\sum_{j=1}^m|z_j|^*\,$. The space $\,X''\,$
provided with the norm (\ref{1.7}) is called the {\sl
Marcinkiewicz\/} space associated with $\,\xb\,$. Since the set
$\,\tilde P_\xb\,$ is $\,\Tf_\wrm(X,X')$-bounded, Theorem
\ref{T1.1} implies that the $\,\Tf_\wrm(X',X)$-bounded set
$\,\{\xb'\in X':\|\xb'\|_L<1\}\,$ absorbs any other
$\,\Tf_\wrm(X',X)$-bounded subset of $\,X'\,$. Therefore, in view
of Remark \ref{R1.3}, the strong topology $\,\Tf_\brm(X'',X')\,$
is generated by the norm $\,\|\cdot\|_M\,$. The Mackey topology
$\,\Tf_\mrm(X'',X')\,$ is strictly coarser than
$\,\Tf_\brm(X'',X')\,$ as $\,\|\yb-\yb^{(m)}\|_M\,$ may be equal
to $\,\|\yb\|_M\,$ for all $\,m\,$.
\end{example}

\begin{remark}\label{R1.11}
Let $\,X_{P_\xb}\,$ be the linear space spanned by $\,P_\xb\,$.
Then $\,X''_{P_\xb}\,$ is the minimal symmetric perfect space
which contains $\,\xb\,$. Obviously,
\begin{enumerate}
\item[(1)]
if $\,\xb\,$ is unbounded then $\,X'_{P_\xb}=\R_{00}^\infty\,$ and
$\,X''_{P_\xb}=\R^\infty\,$;
\item[(2)]
if $\,\xb\in l^\infty\,$ but $\,\xb\not\in\R_0^\infty\,$ then
$\,X'_{P_\xb}=l^1\,$ and $\,X''_{P_\xb}=l^\infty\,$;
\item[(3)]
if $\,\xb\in\R_0^\infty\,$ but $\,\xb\not\in l^1\,$ then
$\,X'_{P_\xb}\,$ is the Lorentz space and $\,X''_{P_\xb}\,$ is the
Marcinkiewicz space associated with $\,\xb\,$ (see Example
\ref{E1.10});
\item[(4)]
if $\,\xb\in l^1\,$ but $\,\xb\not\in\R_{00}^\infty\,$ then
$\,X'_{P_\xb}=l^\infty\,$ and $\,X''_{P_\xb}=l^1\,$;
\item[(5)]
if $\,\xb\in\R_{00}^\infty\,$ then $\,X'_{P_\xb}=\R^\infty\,$ and
$\,X''_{P_\xb}=\R_{00}^\infty\,$.
\end{enumerate}
\end{remark}

\begin{remark}\label{R1.12}
If $\,\xb\not\in\R_{00}^\infty\,$ and $\,P_\xb\subset X\,$ then
$\,X'\subseteq l^\infty\,$ and $\,l^1\subseteq X''\,$. Therefore
for every $\,\Tf_\brm(X,X')$-continuous seminorm $\,p\,$ on
$\,X\,$ there exists a constant $\,C_p\,$ such that $\,p(\xb)\le
C_p\,\|\xb\|_{l^1}\,$ for all $\,\xb\in X\bigcap l^1\,$.
\end{remark}

\section{Birkhoff's theorem}

\subsection{Notation and definitions}
Let $\,\Gb=\{G_1,G_2,\ldots\}\,$ be a family of countable sets
$\,G_k\,$ which may have non-empty intersections. Define a simple
graph $G$ as follows: the set of vertices of $G$ coincides with
$\bigcup_kG_k$ and two vertices are joined by an edge in $G$ if
and only if they belong to the same set $G_k$. Then $G_k$ become
complete subgraphs of $G$. Throughout this section we denote by
$g$ (with or without indices) the vertices of $G$ or, in other
words, the elements of $\bigcup_kG_k$. Let
\begin{enumerate}
\item[]
$\,\WC\,$ be the linear space of real-valued functions $\,\wb\,$
on $G\,$,
\item[]
$\WC_+\,$ be the cone of non-negative functions $\,\wb\in\WC\,$
and
\item[]
$\,\WC_0\,$ be the set of functions $\,\wb\in\WC\,$ which take
only finitely many non-zero values.
\end{enumerate}
We shall call $\,\wb\in\WC_+\,$ {\sl weights\/} over $\,G\,$ and
denote by $\,\wb(g)\,$ the weight assigned to $\,g\in G\,$ (that
is, the value of $\wb$ at $g$). If $\,\wb\in\WC\,$, let
\begin{enumerate}
\item[]
$G_\wb\,$ be the subgraph of $G$ which includes all vertices
$\,g\in G\,$ such that $\,\wb(g)\ne0\,$ and all edges joining
these vertices.
\end{enumerate}

Let $\Gb_1$ be an arbitrary subset of $\Gb$. We shall say that a
weight $\wb\in\WC_+$ is {\sl $\Gb_1$-stochastic\/} if
$\,\sum_{g\in G_k}\wb(g)\le1\,$ for every $\,G_k\in\Gb\,$ and
$\,\sum_{g\in G_k}\wb(g)=1\,$ for every $\,G_k\in\Gb_1\,$. Denote
by $\SC^{\Gb_1}$ the convex set of all $\Gb_1$-stochastic weights
and let $\,\PC^{\Gb_1}\,$ be the set of $\Gb_1$-stochastic weights
taking only the values 0 and 1. Clearly, $\,\wb\in\PC^{\Gb_1}\,$
if and only if the restriction of $\wb$ to every subset $G_k$
takes at most one value 1, all other values being 0, and $\wb$
does take the value 1 at some vertex $g\in G_k$ whenever
$G_k\in\Gb_1$. If $\,\Gb_1\subseteq\Gb'_1\subseteq\Gb\,$ then
$\,\PC^{\Gb'_1}\subseteq\PC^{\Gb_1}\subseteq\SC^{\Gb_1}\,$ and
$\,\SC^{\Gb'_1}\subseteq\SC^{\Gb_1}\,$.

\begin{remark}\label{R2.1}
The weights $\wb\in\SC^\Gb$ and $\wb\in\SC^\emptyset$ are said to
be {\sl stochastic\/} and, respectively, {\sl sub-stochastic\/}. A
stochastic weight $\wb$ can be considered as a family of
probability measures $\wb_{(k)}:=\left.\wb\right|_{G_k}$ on the
sets $G_k$ such that $\wb_{(k)}=\wb_{(j)}$ on $G_k\bigcap G_j$.
\end{remark}

Since the set of vertices is countable, $\,\WC\,$ can be
identified with the sequence space $\,\R^\infty\,$ (or with its
subspace if $\,G\,$ is finite). Further on we use definitions and
notation introduced in Section 1.

\subsection{Extreme points}
We shall say that a path $g_0\to g_1\to\dots\to g_l$ in $G$ is
\begin{enumerate}
\item[]
{\sl admissible\/} if no three adjacent vertices in this path
belong to the same set $\,G_k\in\Gb\,$;
\item[]
a {\sl cycle\/} if $\,g_0=g_l\,$ and the number of distinct
vertices $g_j$ is not smaller than 3 (that is, $\,g_0\to g_1\to
g_2=g_0\,$ is not a cycle).
\end{enumerate}

\begin{proposition}\label{P2.2}
Every two vertices lying in the same connected component of  $G$
can be joined by an admissible path. If there are no admissible
cycles and, in addition,
\begin{enumerate}
\item[(c$_1$)]
the intersection $\,G_k\bigcap G_l\,$ of two distinct sets
$\,G_k,G_l\in\Gb\,$ contains at most one vertex of $\,G\,$
\end{enumerate}
then this admissible path is unique.
\end{proposition}

\begin{proof}
Let $\,g_0\,$ and $\,g_m\,$ belong to the same connected component
of $\,G\,$. Then a path $\,g_0\to g_1\to\dots\to g_m\,$ with the
minimal possible number of vertices is admissible (otherwise we
could obtain a shorter path from $\,g_0\,$ to $\,g_m\,$ replacing
$g_j\to g_{j+1}\to\dots\to g_{j+i}$ with $g_j\to g_{j+i}$). This
proves the first statement.

Let $\,g_1\to g_1\to\dots\to g_m\,$ and $\,g_1\to
g_{n+m}\to\dots\to g_{m+1}\to g_m\,$ be two distinct admissible
paths from $\,g_1\,$ to $\,g_m\,$. Without loss of generality we
may assume that these paths have only two common vertices
$\,g_1\,$ and $\,g_m\,$. Then the vertices
$\,g_1,\ldots,g_{m+n}\,$ are distinct and do not belong to the
same set $\,G_k\,$. Consider the graph $\,\GC\,$ formed by all
these vertices and all joining them edges. Let $\,\tilde
g_1\to\tilde g_2\to\dots\to\tilde g_{l+1}=\tilde g_1\,$ be a cycle
in $\,\GC\,$ with the minimal possible number of vertices which do
not belong to the same set $\,G_k\,$ (since $\,\GC\,$ contains at
least one cycle $\,g_1\to g_2\to\dots\to g_{m+n}\,$ with this
property, such a `minimal' cycle exists). The condition (c$_1$)
implies that this cycle is admissible. Indeed, if two non-adjacent
vertices $\,\tilde g_i\,$ and $\tilde g_{i+j}\,$ in this path are
joined by an edge then all vertices of the cycle $\,\tilde g_i\to
\tilde g_{i+1}\to\dots\to\tilde  g_{i+j}\to\tilde  g_i\,$ belong
to some set $\,G_k\in\Gb\,$ and all vertices of the cycle
$\,\tilde g_{i+j}\to\tilde  g_{i+j+1}\to\dots\to\tilde  g_l\to
\tilde g_1\to\dots\to\tilde  g_i\to\tilde  g_{i+j}\,$ belong to a
distinct set $\,G_l\,$, in which case the intersection
$\,G_k\bigcap G_l\,$ contains at least two elements $\,\tilde
g_i\,$ and $\,\tilde g_{i+j}\,$. This proves the second statement.
\end{proof}

Further on we shall be assuming that
\begin{enumerate}
\item[({\bf g}$_1$)] every vertex of $G$ belongs to at most two
sets $G_k$,
\item[({\bf g}$_2$)] every admissible cycle in $G$ has an even number of
vertices.
\end{enumerate}

If the conditions ({\bf g}$_1$) and ({\bf g}$_2$) are fulfilled
then $\Gb$ can be split into two groups
$\Gb^+=\{G_1^+,G_2^+,\ldots\}$ and $\Gb^-=\{G_1^-,G_2^-,\ldots\}$
in such a way that any two sets from the same group do not have
common elements (two sets $G_k$ and $G_j$ belong to the same group
if every admissible path $\,G_k\ni g_0\to g_1\to\dots\to
g_{l-1}\to g_l\in G_j\,$ in $G$ with $g_1\not\in G_k$ and
$g_{l-1}\not\in G_j$ has an even number of vertices). The
intersection $G_k^+\bigcap G_j^-$ may consist of several elements
or be empty, and every set $G_k^\pm$ may contain a `tail' subset
$\tilde G_k^\pm$ which does not have common elements with any
other set $G_j$.

In view of the following example, all results of this section are
valid for finite and infinite matrices which we shall discuss in
more detail in Section 3.

\begin{example}\label{E2.3}
Let $\Gb$ satisfy ({\bf g}$_1$) and ({\bf g}$_2$) and
$\,\Gb^\pm\,$ be defined as above. Denote by $m_\pm$ the number of
sets $G_k^\pm$ lying in $\Gb^\pm\,$; we allow $\,m_+=\infty\,$ and
(or) $\,m_-=\infty\,$. If every intersection $G_k^+\bigcap G_j^-$
consists of one element and all the tail subsets $\tilde G_k^\pm$
are empty then $\,\WC\,$ is isomorphic to the linear space of
$\,m_+\times m_-$-matrices. Indeed, the value of $\,\wb\in\WC\,$
at the vertex $g\in G_k^+\bigcap G_j^-$ can be considered as the
entry of an $\,m_+\times m_-$-matrix at the intersection of its
$j$th row and $k$th column. In this case $\,\SC^\Gb\,$,
$\,\SC^\emptyset\,$ and $\PC^{\Gb}$ are the sets of doubly
stochastic, sub-stochastic and permutation matrices respectively.
\end{example}

If $G$ is a general family of sets satisfying ({\bf g}$_1$) and
({\bf g}$_2$) then one can think of $\WC$ as a space of
matrices which may have `multiple' or `forbidden' entries and
`tails' $\tilde G_k^\pm$ attached to their rows and columns.

\begin{theorem}\label{T2.4}
Let the conditions {\rm({\bf g}$_1$)} and {\rm({\bf g}$_2$)} be
fulfilled and let $\,\VC\,$ be a normal conic subset of $\WC$.
Then $\,\ex(\SC^{\Gb_1}\bigcap\VC)=\PC^{\Gb_1}\bigcap\VC\,$.
\end{theorem}

\begin{proof}
Obviously,
$\,\PC^{\Gb_1}\bigcap\VC\subset\ex(\SC^{\Gb_1}\bigcap\VC)\,$. In
order to prove the converse, let us consider a weight
$\,\wb\in\SC^{\Gb_1}\bigcap\VC\,$ such that $\,\wb(g')\in(0,1)\,$
for some $\,g'\in G\,$ and show that
$\,\wb\not\in\ex(\SC^{\Gb_1}\bigcap\VC)\,$. Let $\,G'\,$ be the
connected component of $\,G_\wb\,$ containing the vertex $\,g'\,$.
Then $\,\wb(g)\in(0,1)\,$ at every vertex $\,g\in G'\,$.

\medskip\noindent({\bf1}) \
Assume that, for some $\,k\ne l\,$, the intersection $\,G'\bigcap
G_k\bigcap G_l\,$ contains two distinct vertices $\,g_1\,$ and
$\,g_2\,$. Let $\,\wb_\ve^\pm(g_j)=\wb(g_j)\pm(-1)^j\ve\,$ and
$\,\wb_\ve^\pm(g)=\wb(g)\,$ whenever $\,g\ne g_j\,$, $j=1,2\,$.
Then $\,\wb=\frac12(\wb_\ve^++\wb_\ve^-)\,$ and, in view of ({\bf
g}$_1$), $\,\wb_\ve^\pm\in\SC^{\Gb_1}\bigcap\VC\,$ provided that
$\,\ve>0\,$ is sufficiently small. Therefore without loss of
generality we can assume that $\,G'\,$ satisfies (c$_1$).

\medskip\noindent({\bf2}) \
Similarly, if $G'$ contains an admissible cycle $\,\GC=g_0\to
g_2\to\dots\to g_n=g_0\,$, let
$\,\wb_\ve^\pm(g_j)=\wb(g_j)\pm(-1)^j\ve\,$ and
$\,\wb_\ve^\pm(g)=\wb(g)\,$ whenever $\,g\not\in\GC\,$. The
condition ({\bf g}$_2$) implies that $\,\wb_\ve^+\,$ and
$\,\wb_\ve^-\,$ are correctly defined weights over $\,G\,$. We
have $\,\wb=\frac12(\wb_\ve^++\wb_\ve^-)\,$ and
$\,\wb_\ve^{\pm}\in\WC_+\bigcap\VC$ provided that $\ve$ is
sufficiently small. In view of ({\bf g}$_1$), if $g_j\in G_k$ then
one of the adjacent vertices $g_{j-1},g_{j+1}$ belongs to $G_k$
and the other does not. This implies that $\,\sum_{g\in
G_k}\wb_\ve^\pm(g)=\sum_{g\in G_k}\wb(g)$ for every $\,k\,$.
Therefore $\,\wb_\ve^{\pm}\in\SC^{\Gb_1}\bigcap\VC\,$.

\medskip\noindent({\bf3}) \
Finally, let us assume that $G'$ does not contain admissible
cycles and satisfies (c$_1$). Then, by Proposition \ref{P2.2},
every two vertices $g_0,g_l\in G'$ are joined by a unique
admissible path. Let us fix $g_0\in G'$ and denote by $\GC_n$ the
set of vertices in $G'$ obtained from $g_0$ by moving along all
admissible paths with $n$ edges. Then for each $k=1,2,\ldots$
there exists $n\ge 0$ such that
$\,G_k\subseteq\GC_n\bigcup\GC_{n+1}\,$. Moreover, if
$\,G_k\subseteq\GC_n\bigcup\GC_{n+1}\,$ then the intersection
$\,G_k\bigcap\GC_n\,$ consists of one element $\,g_{k,n}\,$.
Indeed, if there are two distinct admissible paths $\,g_0\to
g_1\to\dots\to g_{k,n}\,$ and $\,g_0\to g'_1\to\dots\to
g'_{k,n}\in G_k\bigcap\GC_n\,$ then $\,g_0\,$ and $\,g_{k,n}\,$
can be joined by the two distinct admissible paths $\,g_0\to
g_1\to\dots\to g_{k,n}\,$ and $\,g_0\to g'_1\to\dots\to
g'_{k,n}\to g_{k,n}\,$.

If $\,g_0\in G_k\,$ then $\,G_k\subseteq\GC_0\bigcup\GC_1\,$ and
$\,g_{k,0}=g_0\,$. Let us denote
\begin{equation}\label{2.1}
\ve_{k,0}:=
\min\left\{\frac12\,,\;\frac{1-\wb(g_0)}{2\,\wb(g_0)}\right\}\,,
\qquad\ve_{k,n+1}:=
\frac{\ve_{k,n}\,\wb(g_{k,n})}{1-\wb(g_{k,n})}\,,
\end{equation}
where $\,n=0,1,2,\ldots\,$ and $k$ is such that
$\,G_k\subseteq\GC_n\bigcup\GC_{n+1}\,$. Since
$\,\wb\in\SC^{\Gb_1}\,$, we have
$\,\wb(g_{k,n})+\wb(g_{k,n+1})\le1\,$ and, consequently,
$$
\frac{\wb(g_{k,n})}{1-\wb(g_{k,n})}\ \le\
\frac{1-\wb(g_{k,n+1})}{\wb(g_{k,n+1})}\,.
$$
Using these inequalities, one can easily prove by induction in $n$
that
\begin{equation}\label{2.2}
\ve_{k,n}\ \le\
\min\left\{\frac12\,,\,\frac{1-\wb(g_{k,n})}{2\,\wb(g_{k,n})}\right\}\,.
\end{equation}

Consider two sequences of weights $\wb_{\ve,n}^+$ and
$\wb_{\ve,n}^-$ such that
\begin{enumerate}
\item[]
$\wb_{\ve,0}^\pm(g_0):=(1\pm\ve_{k,0})\,\wb(g_0)\,$ and
$\,\wb_{\ve,0}^\pm(g):=\wb(g)$ for all $g\ne g_0\,$,
\item[]
$\wb_{\ve,n+1}^\pm(g):=\wb_{\ve,n}^\pm(g)\,$ for all
$\,g\in\bigcup_{j\le n}\GC_j\,$,
\item[]
$\wb_{\ve,n+1}^\pm(g):=\wb(g)\,$ whenever $\,g\not\in\bigcup_{j\le
n+1}\GC_j\,$,
\item[]
if $\,G_k\subseteq\GC_n\bigcup\GC_{n+1}\,$ then
$\,\wb_{\ve,n}^\pm(g_{k,n}):=(1\pm\ve_{k,n})\,\wb(g_{k,n})\,$ and
\newline
$\,\wb_{\ve,n+1}^\pm(g):=(1\mp\ve_{k,n+1})\wb(g)\,$ whenever
$\,g\in G_k\bigcap\GC_{n+1}\,$ and $\,g\ne g_{k,n}\,$.
\end{enumerate}
Obviously,
$\,\wb(g)=\frac12(\wb_{\ve,n}^+(g)+\wb_{\ve,n}^-(g))\,$. The
estimates (\ref{2.2}) imply that
$\,\wb_{\ve,n}^\pm\in\WC_+\bigcap\VC\,$. Finally, if
$\,G_k\subseteq\GC_n\bigcup\GC_{n+1}\,$ and $\,\sum_{g\in
G_k}w(g)=t\,$ then
$$
\sum_{g\in G_k}\wb_{\ve,n}^\pm(g) =(1\pm\ve_{k,n})\,\wb(g_{k,n})\
+\ (1\mp\ve_{k,n+1})\,(t-\wb(g_{k,n}))\ =\ t\mp\ve_{k,n+1}(t-1)\,.
$$
This identity and the estimates $\,\ve_{k,n},\ve_{k,n+1}\le1/2\,$
imply that
\begin{equation}\label{2.3}
\frac{t+1}2\;-\;\frac{1-t}{2\,(1-w(g_{k,n}))}\ \le\ \sum_{g\in
G_k}\wb_{\ve,n}^\pm(g)\ \le\ \frac{t+1}2\,.
\end{equation}
Let $\wb_\ve^\pm(g):=\lim_{n\to\infty}\wb_{\ve,n}^\pm(g)$. Then
$\wb_0=\frac12(\wb_\ve^++\wb_\ve^-)$ and, in view of (\ref{2.2})
and (\ref{2.3}), $\wb_\ve^\pm(g)\in\SC^{\Gb_1}\bigcap\VC\,$.

Thus, under conditions of the theorem, a weight
$\,\wb\in(\SC^{\Gb_1}\setminus\PC^{\Gb_1})\bigcap\VC\,$ can always
be represented as a convex combination of two other weights from
$\,\SC^{\Gb_1}\bigcap\VC\,$ and therefore is not an extreme point.
\end{proof}

\begin{remark}\label{R2.5}
If the condition ({\bf g}$_2$) is not fulfilled then an extreme
point of $\,\SC^{\Gb_1}\,$ does not necessarily belong to
$\,\PC^{\Gb_1}\,$. The simplest example is $\,G_1=\{g_1,g_2\}\,$,
$\,G_2=\{g_2,g_3\}\,$, $\,G_3=\{g_3,g_1\}\,$ and
$\,\Gb=\{G_1,G_2,G_3\}\,$. In this case $\SC^\Gb$ consists of one
weight which takes the value $\frac12$ at each vertex.
\end{remark}

\begin{remark}\label{R2.6}
The sets $\SC^{\Gb_1}$ and $\PC^{\Gb_1}$ may well be very
poor or even empty.
However, even in this situation Theorem \ref{T2.4} may be useful.
In particular, by the Krein--Milman theorem, under conditions of
Theorem \ref{T2.4} we have
$$
\SC^{\Gb_1}\bigcap\WC_0\ =\ \ch \PC^{\Gb_1}\bigcap\WC_0\,.
$$
Therefore $\SC^{\Gb_1}\bigcap\WC_0=\emptyset$
whenever $\PC^{\Gb_1}\bigcap\WC_0=\emptyset$.
\end{remark}

\begin{remark}\label{R2.7}
If the conditions ({\bf g}$_1$) and ({\bf g}$_2$) are fulfilled
and $\,\VC\,$ is a normal linear subspace of $\,\WC\,$ then every
extreme point
$\,\wb\in\ex(\SC^{\Gb_1}\bigcap\VC)=\PC^{\Gb_1}\bigcap\VC\,$ is
$\,\Tf_\mrm(\VC,\VC')$-exposed. Indeed, if $\,\wb'(g)>0\,$
whenever $\,\wb(g)=1\,$, $\,\wb'(g)<0\,$ whenever $\,\wb(g)=0\,$
and $\,\wb'\in\VC'\,$ then we have
$\,\langle\wb,\wb'\rangle>\langle\tilde\wb,\wb'\rangle\,$ for all
$\,\tilde\wb\in\SC^{\Gb_1}\bigcap\VC\,$.
\end{remark}

\subsection{Topologies on the space of stochastic weights}
The aim of this subsection is to describe locally convex
topologies $\,\Tf\,$ on a linear subspace
$\,\VC\supset\PC^{\Gb_1}\,$ such that the $\,\Tf$-closure of
$\,\ch\PC^{\Gb_1}\,$ coincides with $\,\SC^{\Gb_1}\bigcap\VC\,$.
By Fatou's lemma we always have
$\,\cch\PC^\emptyset\subset\SC^\emptyset\,$ (as $\,\Tf\,$ is finer
than $\,\Tf_0\,$). Tychonoff's theorem and Fatou's lemma also
imply that the set $\,\SC^\emptyset\,$ is $\,\Tf_0$-compact.
Therefore, in view of Theorem \ref{T2.4} and the Krein--Milman
theorem, under the conditions ({\bf g}$_1$) and ({\bf g}$_2$) we
have $\,\SC^\emptyset=\cch\PC^\emptyset\,$, where the closure is
taken in the topology of element-wise convergence $\,\Tf_0\,$.
However, if $\,\Gb_1\,$ contains an infinite set $\,G_k\,$ then
the set $\,\SC^{\Gb_1}\,$ is not $\,\Tf_0$-closed and, by Theorem
\ref{T1.2}, is not $\,\Tf$-compact whenever the functional
$\,\wb\to\sum_{g\in G_k}\wb(g)\,$ is $\,\Tf$-continuous. In this
case {\bf(ii)} does not directly follow from {\bf(i)} and the
Krein--Milman theorem.

\begin{definition}\label{D2.8}
Denote by $\,\VC_\PC\,$ and $\,\VC_\SC\,$ the normal covers of the
subspaces spanned by $\,\PC^\emptyset\,$ and $\,\SC^\emptyset\,$
respectively. If $\,\wb\in\VC_\SC\,$, let $\,\wb_{(k)}\,$ be the
restriction of $\,\wb\,$ to $\,G_k\,$ and
$\,p_k(\wb):=\|\wb_{(k)}\|_{l^1}\,$.
\end{definition}

\begin{lemma}\label{L2.9}
Let us enumerate the sets $\,G_k\,$ in an arbitrary way and define
$\,F_n:=\bigcup_{k=1}^nG_k\,$. If $\,D'\,$ is a
$\,\Tf_\wrm(\VC'_\PC,\VC_\PC)$-compact subset of $\,\VC'_\PC\,$
then
\begin{equation}\label{2.4}
\sup_{\wb\in\PC^{\emptyset},\,\wb'\in D'} \,\sum_{g\in G\setminus
F_n}|\wb(g)\,\wb'(g)|\ \to\ 0\,,\qquad n\to\infty\,,
\end{equation}
whenever $\,\Gb\,$ satisfies {\rm({\bf g}$_1$)} and
\begin{equation}\label{2.5}
\sup_{\wb\in\SC^{\emptyset},\,\wb'\in D'} \,\sum_{g\in G\setminus
F_n}|\wb(g)\,\wb'(g)|\ \to\ 0\,,\qquad n\to\infty\,,
\end{equation}
whenever $\,\Gb\,$ satisfies {\rm({\bf g}$_1$)} and {\rm({\bf
g}$_2$)}.
\end{lemma}

\begin{proof}
If the conditions ({\bf g}$_1$) and ({\bf g}$_2$) are fulfilled then
$\,\SC^{\emptyset}\,$ coincides with the $\,\Tf_0$-closure of
$\,\ch\PC^{\emptyset}\,$. Therefore, in view of Fatou's lemma, it
is sufficient to prove only the first statement.

If (\ref{2.4}) is not true then there exists $\,\de>0\,$ and two
sequences of weights $\,\{\wb_n\}\subset\PC^{\emptyset}\,$ and
$\,\{\wb'_n\}\subset D'\,$ such that $\,\sum_{g\in G\setminus
F_n}|\wb_n(g)\,\wb'_n(g)|\ge\de>0\,$ for all $n=1,2,\ldots$ Let
$n_\star$ be the minimal positive integer satisfying the estimate
$\,\sum_{g\in F_{n_\star}\setminus F_n}
|\wb_n(g)\,\wb'_n(g)|\ge\de/2\,$ and
$\wb^\star_n\in\WC_0\bigcap\PC^\emptyset$ be the weight which
takes the same values as $\wb_n$ on $\,F_{n_\star}\setminus F_n\,$
and vanishes outside $\,F_{n_\star}\setminus F_n\,$. In view of
({\bf g}$_1$), there exists a positive integer $n^\star>n_\star$
such that $\left.\wb^\star_n\right|_{G_k}\equiv0$ for all $k\ge
n^\star$. Let us take an arbitrary $n_1$ and define
$n_{j+1}:=n_j^\star\,$, where $j=1,2,\ldots$ Then for each $\,g\in
G\,$ the sum $\,\wb^\star(g):=\sum_j\wb_{n_j}^\star(g)\,$ is equal
either to 0 or to 1 and $\,\sum_{g\in G_k}\wb^\star(g)\le 1\,$,
$\,\forall k=1,2,\ldots$ Therefore the corresponding weight
$\,\wb^\star\,$ belongs to $\,\PC^{\emptyset}\,$. On the other
hand, $\,n_j\to\infty\,$ and
$$
\sum_{g\in F_{n_j^\star}\setminus F_{n_j}}
|\wb^\star(g)\,\wb'_{n_j}(g)|\ =\sum_{g\in F_{n_j^\star}\setminus
F_{n_j}} |\wb_{n_j}(g)\,\wb'_{n_j}(g)|\ \ge\ \de/2\,,
$$
which contradicts to Theorem \ref{T1.2}.
\end{proof}

We do not assume in Lemma \ref{L2.9} that $\,D'\subset\VC_\SC\,$.
Therefore, for each fixed $\,n\,$, the supremum in (\ref{2.5}) may
well be $\,+\infty\,$. However, under conditions ({\bf g}$_1$) and
({\bf g}$_2$), it eventually becomes finite and converges to zero
as $\,n\to\infty\,$.

\begin{lemma}\label{L2.10}
If the condition {\rm({\bf g}$_1$)} is fulfilled and $\,D'\,$ is a
$\,\Tf_\wrm(\VC'_\SC,\VC_\SC)$-bounded
$\,\Tf_\wrm(\VC'_\PC,\VC_\PC)$-compact subset of $\,\VC'_\SC\,$
then the weights $\,\wb'\in D'\,$ are uniformly bounded.
\end{lemma}

\begin{proof}
Let $\,F_n\,$ be defined as in Lemma \ref{L2.9}. If the
restrictions of weights $\,\wb'\in D'\,$ to $\,F_n\,$ are not
uniformly bounded then, for some $\,k\le n\,$, their restrictions
to $\,G_k\,$ form an unbounded subset of $\,l^\infty\,$. This
implies that the set $\,D'\,$ is not
$\,\Tf_\wrm(\VC'_\SC,\VC_\SC)$-bounded.

Assume that there exist sequences $\,\{g_j\}_{j=1,2,\ldots}\in
G\,$ and $\,\{\wb'_j\}_{j=1,2,\ldots}\in D'\,$ such that
$\,\wb'_j(g_j)\to\infty\,$ as $\,j\to\infty\,$ and
$\,\{g_j\}\not\subset F_n\,$ for any finite $\,n\,$. Since ({\bf
g}$_1$) holds true, every vertex $\,g\,$ belongs only to finitely
many sets $\,G_k\,$ and we can find a subsequence
$\,\{g_{j_i}\}_{i=1,2,\ldots}\,$ with at most one entry at each
set $\,G_k\,$. If $\,\wb(g_{j_i})=1\,$ and $\,\wb(g)=0\,$ whenever
$\,g\not\in\{g_{j_i}\}\,$ then $\,\wb\in\PC^\emptyset\,$ and
$\,\sum_i|\wb(g_{j_i})\,\wb'_{j_i}(g_{j_i})|=\infty\,$. Therefore,
by Theorem \ref{T1.2}, the set $\,D'\,$ is not
$\,\Tf_\wrm(\VC'_\PC,\VC_\PC)$-compact.
\end{proof}

\begin{corollary}\label{C2.11}
If the conditions {\rm({\bf g}$_1$)} and {\rm({\bf g}$_2$)} are
fulfilled then the strong topology
$\,\Tf_\brm(\VC_\SC,\VC'_\SC)\,$ is generated by the norm
\begin{equation}\label{2.6}
\|\wb\|_\SC\ :=\ \sup_k p_k(\wb)
\end{equation}
\end{corollary}

\begin{proof}
Since the norm (\ref{2.6}) is lower $\,\Tf_0$-semicontinuous, it
is $\,\Tf_\brm(\VC_\SC,\VC'_\SC)$-continuous (see Remark
\ref{R1.5}). The set
$$
\SC^\emptyset\ =\ \{\wb\in\VC_\SC:\|\wb\|_\SC<1\}
$$
is absorbing and, in view of (\ref{2.5}) and Lemma \ref{L2.10}, is
$\,\Tf_\mrm(\VC_\SC,\VC'_\SC)$-bounded. By Theorem \ref{T1.1}, this
set is $\,\Tf_\brm(\VC_\SC,\VC'_\SC)$-bounded, which implies that
every $\,\Tf_\brm(\VC_\SC,\VC'_\SC)$-continuous seminorm is
continuous in the norm topology.
\end{proof}

\begin{corollary}\label{C2.12}
Let $\,\Tf\,$ be the locally convex topology on $\,\VC_\SC\,$
generated by the seminorms $\,p_k\,$, $\,k=1,2,\ldots\,$ If the
conditions {\rm({\bf g}$_1$)} and {\rm({\bf g}$_2$)} are fulfilled
then the Mackey topology $\,\Tf_\mrm(\VC_\SC,\VC'_\SC)\,$ is finer
than $\,\Tf\,$ and coincides with $\,\Tf\,$ on every
$\,\Tf_\brm(\VC_\SC,\VC'_\SC)$-bounded subset of $\,\VC_\SC\,$.
\end{corollary}

\begin{proof}
The seminorms $\,p_k\,$ are lower
$\,\Tf_0$-semicontinuous and satisfy (\ref{1.3}). Therefore, by
Remark \ref{R1.5}, $\,\Tf_\mrm(\VC_\SC,\VC'_\SC)\,$ is finer than
$\,\Tf\,$. On the other hand, if $\,\Om\,$ is a bounded subset of
$\,\VC_\SC\,$ then, in view of (\ref{2.5}) and Lemma \ref{L2.10},
for every Mackey seminorm $\,p\,$ on $\,\VC_\SC\,$, every
$\,\xb\in\Om\,$ and every $\,\ve>0\,$ there exist a positive
integer $\,m\,$ and $\,\de>0\,$ such that
$$
\{\yb\in\Om:p_k(\xb-\yb)<\de\,,\forall k=1,2,\ldots,m\}\
\subseteq\ \{\yb\in\Om:p(\xb-\yb)<\ve\}\,.
$$
This implies that every
$\,\Tf_\mrm(\VC_\SC,\VC'_\SC)$-neighbourhood of $\,\xb\,$ in
$\,\Om\,$ contains a $\,\Tf$-neighbourhood.
\end{proof}

\begin{remark}\label{R2.13}
If the conditions ({\bf g}$_1$), ({\bf g}$_2$) are fulfilled and
$\,G\,$ does not coincide with the union of a finite collection of
the sets $\,G_k\,$ then the topology $\,\Tf\,$ generated by the
seminorms $\,p_k\,$ is strictly coarser than
$\,\Tf_\mrm(\VC_\SC,\VC'_\SC)\,$. Indeed, in this case there
exists a sequence of weights $\,\wb_n\in\VC_\PC\,$ such that
$\,p_k(\wb_n)=0\,$ for all $\,k<n\,$ and $\,p_n(\wb_n)\to\infty\,$
as $\,n\to\infty\,$. This sequences converges to the zero weight
in the topology $\,\Tf\,$ but is not
$\,\Tf_\brm(\VC_\SC,\VC'_\SC)\,$-bounded and, consequently, is not
$\,\Tf_\mrm(\VC_\SC,\VC'_\SC)\,$-convergent.
\end{remark}

In the rest of this section we shall be assuming that
\begin{enumerate}
\item[({\bf g}$_3$)]
one can enumerate the sets $\,G_j\,$ in such a way that either
$\,G=F_n\,$ or $\,G_{n+1}\not\subset F_n\,$ for all sufficiently
large $\,n\,$, where $\,F_n:=\bigcup_{k\le n}G_k\,$.
\end{enumerate}
Every finite collection $\,\Gb=\{G_1,G_2,\ldots,G_n\}\,$ satisfies
({\bf g}$_3$). More generally, the condition ({\bf g}$_3$) is
fulfilled whenever the number of finite sets $\,G_k\,$ is finite
and the intersections of every two sets $\,G_j,G_k\in\Gb\,$ is
finite. In particular, ({\bf g}$_3$) is fulfilled for finite and
infinite matrices (see Example \ref{E2.3}).

\begin{lemma}\label{L2.14}
Let the conditions {\rm({\bf g}$_1$)} and {\rm({\bf g}$_3$)} be
fulfilled, $\,G_k\,$ be enumerated as in {\rm({\bf g}$_3$)},
$\,F_n:=\bigcup_{k\le n}G_k\,$ and $\,\Gb_{1,n}\,$ be the
collection of all sets $\,G_k\in\Gb_1\,$ with $\,k\le n\,$. Then
there exists a positive integer $\,n_0\,$ such that for every
$\,n\ge n_0\,$ and every weight $\,\wb\in\PC^{\emptyset}\,$
satisfying
\begin{equation}\label{2.7}
\sum_{g\in G_k}\wb(g)=1\,,\qquad\forall G_k\in\Gb_{1,n}\,,
\end{equation}
one can find a weight $\,\tilde\wb\in\PC^{\Gb_1}\,$ whose
restriction to $\,F_n\,$ coincides with
$\,\left.\wb\right|_{F_n}\,$.
\end{lemma}

\begin{proof}
If for some positive integer $\,n_1\,$ there are no weights
$\,\wb\in\PC^{\emptyset}\,$ satisfying (\ref{2.7}) with
$\,n=n_1\,$ then the lemma automatically holds true for
$\,n_0=n_1\,$. Therefore we can assume without loss of generality
that for each $\,n=1,2,\ldots\,$ there exists a weight
$\,\wb_n\in\PC^{\emptyset}\,$ satisfying (\ref{2.7}).

If $\,G=F_n\,$ for all $\,n\ge n_1\,$ then, in view of ({\bf
g}$_1$), $\,\sum_k\sum_{g\in G_k}\wb_n(g)\le 2n_1\,$. This
estimate and (\ref{2.7}) imply that the set $\,\Gb_{1,n}\,$
contains at most $\,2n_1\,$ elements for each $\,n=1,2,\ldots\,$.
Therefore there exists a positive integer $\,n_0\,$ such that
$\,\Gb_{1,n}=\Gb_{1,n_0}\,$ for all $\,n\ge n_0\,$. In this case
the inclusion $\,\wb\in\PC^\emptyset\,$ and (\ref{2.7}) with
$\,n\ge n_0\,$ imply that $\,\wb\in\PC^{\Gb_1}\,$.

If $\,G_{n+1}\not\subset F_n\,$ for all $\,n\ge n_1\,$ then we
take $\,n_0=n_1\,$. Given $\,n\ge n_0\,$ and a weight
$\,\wb\in\PC^{\emptyset}\,$ satisfying (\ref{2.7}), we choose
 arbitrary vertices $\,g_{n+j}\in G_{n+j}\setminus
F_{n+j-1}\,$ and define $\,\tilde\wb\,$ as follows:
$\,\tilde\wb(g):=\wb(g)\,$ whenever $\,g\in F_n\,$,
$\,\tilde\wb(g_{n+j}):=1\,$ for all $\,j=1,2,\ldots\,$ and
$\,\tilde\wb(g):=0\,$ otherwise. Then $\,\wb=\tilde\wb\,$ on
$\,F_n\,$ and $\,\tilde\wb\in\PC^{\Gb_1}\,$ because $\,\sum_{g\in
G_k}\tilde\wb(g)=1\,$ for all $\,k>n\,$.
\end{proof}

\begin{theorem}\label{T2.15}
Let the conditions {\rm ({\bf g}$_1$)--({\bf g}$_3$)} be fulfilled
and $\,\VC\,$ be a normal subspace of $\,\WC\,$ such that
$\,\VC_\PC\subseteq\VC\subseteq\VC_\SC\,$. Then
\begin{equation}\label{2.8}
\SC^{\Gb_1}\bigcap\VC\ =\ \cch\PC^{\Gb_1}\,,
\qquad\forall\Gb_1\subseteq\Gb\,,
\end{equation}
where the closure is taken in the Mackey topology
$\,\Tf_\mrm(\VC,\VC')\,$.
\end{theorem}

\begin{proof}
Since $\,\VC\subseteq\VC_\SC\,$, the functionals
$\,\wb\to\sum_{g\in G_k}\wb(g)\,$ are
$\,\Tf_\wrm(\VC,\VC')$-continuous and, consequently,
$\,\cch\PC^{\Gb_1}\subseteq \SC^{\Gb_1}\bigcap\VC\,$. If
$\,\wb\not\in\cch\PC^{\Gb_1}\,$ then, by the separation theorem
(see, for example, \cite{K}, Section 20.7), there exist
$\wb'\in\VC'$ and $\,\ve>0\,$ such that
$\,\langle\wb,\wb'\rangle-\langle\tilde\wb,\wb'\rangle>\ve\,$ for
all $\,\tilde\wb\in\ch\PC^{\Gb_1}\,$. Therefore, in order to prove
(\ref{2.8}), it is sufficient to show that for each fixed
$\,\wb\in\SC^{\Gb_1}\bigcap\VC\,$, $\wb'\in\VC'$ and $\ve>0$ one
can find $\,\tilde\wb\in\ch\PC^{\Gb_1}\,$ such that
$\,\langle\wb,\wb'\rangle-\langle\tilde\wb,\wb'\rangle\le\ve\,$.

Assume that the intersection $\,G_k\bigcap G_l\bigcap G_\wb\,$
contains more than one vertex so that $\,G_k\bigcap G_l\bigcap
G_\wb=\{g_1,g_2,\ldots\}\,$. Since
$\,\sum_j|\wb(g_j)\,\wb'(g_j)|\le\infty\,$ and
$\,\sum_j\wb(g_j)\le1\,$, we have
$\,\wb'(g_i)\ge\sum_j\wb(g_j)\,\wb'(g_j)\,$ for some positive
integer $\,i\,$. If $\,\wb^\star(g):=\wb(g)\,$ whenever
$\,g\not\in G_k\bigcap G_l\bigcap G_\wb\,$,
$\,\wb^\star(g_i):=\sum_j\wb(g_j)\,$ and $\,\wb^\star(g):=0\,$
whenever $\,g\in G_k\bigcap G_l\bigcap G_\wb\,$ but $\,g\ne g_i\,$
then $\,\wb^\star\in\SC^{\Gb_1}\bigcap\VC\,$ and
$\,\langle\wb^\star,\wb'\rangle\ge\langle\wb,\wb'\rangle\,$.
Therefore we can assume without loss of generality that
$\,G_\wb\,$ satisfies the condition (c$_1$) of Proposition
\ref{P2.2}.

Let us enumerate the sets $\,G_k\,$ and define $F_n$ and $\,n_0\,$
as in Lemma \ref{L2.14}. Let $\,n\ge n_0\,$ and $\,G_k\bigcap
G_\wb=\{g_1^k,g_2^k,...\}\,$, where $\,k=1,2,\ldots,n\,$. By ({\bf
g}$_1$), for every $\,g_j^k\,$ there exists at most one positive
integer $\,l\ne k\,$ such that $\,g_j^k\in G_l\,$. Denote
$$
v_n(g_j^k)\ :=\ \begin{cases}\wb(g_j^k)\,,&\text{if}\
g_j^k\not\in\bigcup_{l=n+1}^\infty G_l\,,\\
\sum_{g\in G_l\bigcap F_n}\wb(g)\,,&\text{if}\ g_j^k\in G_l\
\text{for some}\ l>n\,.
\end{cases}
$$
In view of ({\bf g}$_1$) and (c$_1$), we have
$\,\sum_jv_n(g_j^k)\le n\,$. Therefore $\,v_n(g_j^k)\to0\,$ as
$\,j\to\infty\,$ whenever the set $\,G_k\bigcap G_\wb\,$ is
infinite. If $\,G_k\bigcap G_\wb\,$ is finite, denote by $\,j_k\,$
the number of elements of $\,G_k\bigcap G_\wb\,$. If $\,G_k\bigcap
G_\wb\,$ is infinite, denote by $\,j_k\,$ the minimal positive
integer such that
\begin{equation}\label{2.9}
v_n(g_{j_k}^k)+\sum_{j>j_k}\wb(g_j^k)\ \le\
1\qquad\text{and}\qquad g_j^k\not\in F_n\setminus
G_k\,,\quad\forall j\ge j_k
\end{equation}
(since $\,v_n(g_j^k)\to0\,$, $\,\sum_j\wb(g_j^k)\le1\,$ and
$\,G_\wb\,$ satisfies (c$_1$), such a minimal integer exists).

Let $\,\wb_n(g):=0\,$ whenever $\,g\not\in F_n\,$ and
$\,\wb_n(g):=\wb(g)\,$ for all $\,g\in F_n\,$. Then
$\,\langle\wb-\wb_n,\wb'\rangle\to0\,$ as $\,n\to\infty\,$ because
the series $\,\sum_{g\in G}\wb(g)\,\wb'(g)\,$ is absolutely
convergent. Let $\,m\ge\max\{j_1,j_2,\ldots,j_n\}\,$,
$\,\wb_{n,m}(g):=0\,$ whenever $\,\wb_n(g)=0\,$ and
$$
\wb_{n,m}(g_j^k)\ :=\ \begin{cases}0
&\text{if $j>m$,}\\
\wb(g_{j_k}^k)+\sum_{j>m}\wb(g_j^k)
&\text{if $j=j_k\le m$,}\\
\wb(g)&\text{if $j\le m$ and $j\ne j_k\,$.}
\end{cases}
$$
Then $\,\langle\wb_n-\wb_{n,m},\wb'\rangle\to0\,$ as
$\,m\to\infty\,$ for each fixed $n$ because the series
$\,\sum_j\wb(g_j^k)\,\wb'(g_j^k)\,$ are absolutely convergent and
$\,\sum_{j>m}\wb(g_j^k)\to0\,$.

The weight $\,\wb_{n,m}\,$ vanishes outside a finite subset of
$\,G\,$ and, in view of (\ref{2.9}), belongs to
$\,\SC^\emptyset\,$ and satisfies the condition (\ref{2.7}).
Applying Theorem \ref{T2.4} to the family of sets $\,\{G_1\bigcap
G_{\wb_{n,m}},\ldots,G_n\bigcap G_{\wb_{n,m}}\}\,$ and then the
Krein--Milman theorem, we see that $\,\wb_{n,m}\,$ can be
represented as a finite convex combination
$\,\sum_i\al_i\,\wb_{n,m}^{(i)}\,$ of some weights
$\,\wb_{n,m}^{(i)}\in\PC^\emptyset\,$. Obviously, each weight
$\,\wb_{n,m}^{(i)}\,$ also satisfies (\ref{2.7}). By Lemma
\ref{L2.14}, we can find $\,\tilde\wb_{n,m}^{(i)}\in\PC^{\Gb_1}\,$
such that $\,\tilde\wb_{n,m}^{(i)}= \tilde\wb_{n,m}^{(i)}\,$ on
the set $\,F_n\,$. If
$\,\tilde\wb_{n,m}:=\sum_i\al_i\,\tilde\wb_{n,m}^{(i)}\,$ then
$\,\tilde\wb_{n,m}\in\ch\PC^{\Gb_1}\,$ and, in view of
(\ref{2.4}), we have
$\,\langle\wb_{n,m}-\tilde\wb_{n,m},\wb'\rangle<\ve/3\,$ for all
$\,m\ge\max\{j_1,j_2,\ldots,j_n\}\,$ provided that $\,n\,$ is
sufficiently large. Therefore, choosing a sufficiently large
$\,n\ge n_0\,$ and then a sufficiently large
$\,m\ge\max\{j_1,j_2,\ldots,j_n\}\,$, we can make the right hand
side of the identity
$$
\langle\wb,\wb'\rangle-\langle\tilde\wb_{n,m},\wb'\rangle\ =\
\langle\wb-\wb_n,\wb'\rangle+\langle\wb_n-
\wb_{n,m},\wb'\rangle+\langle\wb_{n,m}-\tilde\wb_{n,m},\wb'\rangle
$$
smaller than $\,\ve\,$.
\end{proof}

\begin{remark}\label{R2.16}
Theorem \ref{T2.15} implies that $\,\PC^{\Gb_1}\ne\emptyset\,$
whenever $\,\SC^{\Gb_1}\ne\emptyset$ and $\,\Gb\,$ satisfies ({\bf
g}$_1$)--({\bf g}$_3$). If
$\,G_{n+1}\not\subset\bigcup_{k=1}^nG_j\,$ for all
$\,n=1,2,\ldots\,$ then, using the same procedure as in the proof
of Lemma \ref{L2.14}, one can show that $\,\PC^\Gb\ne\emptyset\,$.
\end{remark}

\begin{remark}\label{R2.17}
If $\,\VC\,$ is a proper normal subspace of $\,\VC_1\,$ then
$\,\VC'_1\,$ is a proper subspace of $\,\VC'\,$ and the Mackey
topology $\,\Tf_\mrm(\VC,\VC')\,$ is strictly finer than
$\,\Tf_\mrm(\VC_1,\VC'_1)\,$. Therefore choosing a smaller space
$\,\VC\,$ in Theorem \ref{T2.15} we obtain a stronger result which
is valid for a narrower class of $\,\Gb_1$-stochastic weights.
\end{remark}

\begin{remark}\label{R2.18}
Taking $\,\VC=\VC_\SC\,$ in Theorem \ref{T2.15} and applying
Corollary \ref{C2.12}, we obtain
$\,\SC^{\Gb_1}=\cch\PC^{\Gb_1}\,$, where the closure is taken in
the topology generated by the seminorms $\,p_k\,$. This topology
is metrizable. Therefore, under conditions ({\bf g}$_1$)--({\bf
g}$_3$), for every $\,\wb\in\SC^{\Gb_1}\,$ there exists a sequence
of weights $\,\wb_n\in\ch\PC^{\Gb_1}\,$ such that
$\,p_k(\wb-\wb_n)\to0\,$ as $\,n\to\infty\,$ for all
$\,k=1,2,\ldots$
\end{remark}

The following simple example shows that, generally speaking,
$\,\cch\PC^\emptyset\,$ does not contain $\,\SC^\Gb\bigcap\VC\,$
if we take the closure in the strong topology
$\,\Tf_\brm(\VC,\VC')$.

\begin{example}\label{E2.19}
Let $\,\Gb\,$ be an infinite collection of mutually disjoint sets
$\,G_k\,$ such that $\,G_k\,$ contains $\,k\,$ elements. Then the
weight $\,\wb\,$ which takes the values $\,k^{-1}\,$ on $\,G_k\,$
belongs to $\,\SC^\Gb\,$. On the other hand, for every weight
$\,\tilde\wb\in\ch\PC^\emptyset\,$ there exists a positive integer
$\,n\,$ such that the number of nonzero entries in
$\,\left.\tilde\wb\right|_{G_k}\,$ does not exceed $\,n\,$ for
every $\,k\,$. Therefore $\,\|\wb-\tilde\wb\|_\SC=1\,$ for all
$\,\tilde\wb\in\ch\PC^\emptyset\,$, where $\,\|\cdot\|_\SC\,$ is
defined by (\ref{2.6}).
\end{example}

The strong closure of the convex hull of the set of permutation
matrices is also strictly smaller then the set of doubly
stochastic matrices \cite{Is}.

\section{Operators generated by stochastic matrices}

\subsection{Notation and definitions}
In the rest of the paper (with the exception of the proof of
Theorem \ref{T3.15}) we shall be assuming that $\,\WC\,$ is the
space of real matrices $\,\wb=\{w_{ij}\}_{i,j=1,2,\ldots}\,$ and
the sets $G_k$ are the rows and columns (see Example \ref{E2.3}).
Recall that in this case $\,\Gb\,$ satisfies the conditions ({\bf
g}$_1$)--({\bf g}$_3$), $\,\SC^\Gb\,$ and $\,\SC^\emptyset\,$ are
the sets of doubly stochastic and sub-stochastic matrices
respectively, $\,\,\PC^\Gb\,$ is the set of permutation matrices
and $\,\PC^\emptyset\,$ is the set of sub-stochastic matrices
whose entries are equal either to 0 or to 1. For the sake of
definiteness we shall consider only infinite matrices; the
corresponding results for finite matrices are much simpler and can
be proved in a similar manner.

Every matrix $\,\wb\in\WC\,$ generates the linear operator
\begin{equation}\label{3.1}
\R^\infty\ni\xb\ \to\ \{\,\sum_{j=1}^\infty
w_{1j}\,x_j\,, \,\sum_{j=1}^\infty w_{2j}\,x_j\,,\ldots\}\ \in\
\R^\infty
\end{equation}
with domain $\,\DC(\wb)\:=
\{\xb\in\R^\infty\,:\sum_{j=1}^\infty|\,\,w_{ij}\,x_j\,|<\infty\,,\
\forall i=1,2,\ldots\}\,$. We shall denote this operator by the
same letter $\,\wb\,$. Obviously, $\,l^\infty\subseteq \DC(\wb)\,$
for all $\,\wb\in\SC^{\emptyset}\,$ and $\,\DC(\wb)=\R^\infty\,$
for all $\,\wb\in\PC^{\emptyset}\,$, but
$\,\DC(\wb)\ne\R^\infty\,$ whenever $\,\wb\,$ has a row with
infinitely many nonzero entries.

\begin{lemma}\label{L3.1}
If $\,X\subseteq l^\infty\,$ is a symmetric perfect space and
$\,X\ne\R_{00}^\infty\,$ then the operator generated by a matrix
$\,\wb\in\SC^\emptyset\,$ maps $\,X\,$ into $\,X\,$ and is
continuous in the topologies $\,\Tf_\wrm(X,X')\,$,
$\,\Tf_\mrm(X,X')\,$ and $\,\Tf_\brm(X,X')\,$.
\end{lemma}

\begin{proof}
Since $\,X\ne\R_{00}^\infty\,$, by Remark \ref{R1.11} we have
$\,X'\subseteq l^\infty\,$. The inclusions $\,X\subseteq
l^\infty\,$, $\,X'\subseteq l^\infty\,$ and (\ref{2.5}) imply that
$\,\xb'\otimes\xb\in\VC'_\SC\,$ for all $\,\xb\in X\,$ and
$\,\xb'\in X'\,$, which means that $\,\wb\,$ maps the perfect
space $\,X\,$ into itself. Similarly, the transposed operator
$\,\wb^T\,$ maps the perfect space $\,X'\,$ into itself. Therefore
$\,|\langle\wb\xb,\xb'\rangle|=|\langle\xb,\wb^T\xb'\rangle|\,$ is
a $\,\Tf_\wrm(X,X')$-continuous seminorm on $\,X\,$ for each
$\,\xb'\in X'\,$ and is a $\,\Tf_\wrm(X',X)$-continuous seminorm
on $\,X'\,$ for each $\,\xb\in X\,$. This implies that $\,\wb\,$
is $\,\Tf_\wrm(X,X')$-continuous and $\,\wb^T\,$ is
$\,\Tf_\wrm(X',X)$-continuous. Since the continuous operator
$\,\wb^T\,$ maps compact sets into compact sets and bounded sets
into bounded sets, the operator $\,\wb\,$ is
$\,\Tf_\mrm(X,X')$-continuous and $\,\Tf_\brm(X,X')$-continuous.
\end{proof}

\begin{remark}\label{R3.2}
Let $\,X\subseteq l^\infty\,$ be a symmetric space,
$\,X\not\subset\R_{00}^\infty\,$ and $\,\{\xb'\otimes\xb\}\,$ be
the set which contains one element $\,\xb'\otimes\xb\,$, where
$\,\xb\in X\,$ and $\,\xb'\in X'\,$. Applying (\ref{2.4}) to
$\,D'=\{\xb'\otimes\xb\}\,$, we see that the set
$\,P_\xb^\emptyset\,$ is $\,\Tf_\wrm(X,X')$-bounded. Therefore
Theorem \ref{T1.1} implies the first statement of Corollary
\ref{C1.7}. If $\,\Om'\,$ is a $\,\Tf_\wrm(X',X)$-compact subset of
$\,X'\,$ then, by Theorem \ref{T1.2}, the set
$\,D':=\bigcup_{\xb'\in\Om'}\{\xb'\otimes\xb\}\,$ is
$\,\Tf_\wrm(\VC'_\PC,\VC_\PC)$-compact. Therefore Theorem
\ref{T1.2} and (\ref{2.4}) imply the second statement of Corollary
\ref{C1.7}.
\end{remark}

\begin{definition}\label{D3.3}
Let $\,\Gb_r\,$ be the set of all rows, $\,\SC^r:=\SC^{\Gb_r}\,$
and $\,\UC\SC^r\,$ be the set of matrices
$\,\wb=\{w_{ij}\}_{i,j=1,2,\ldots}\in\SC^r\,$ such that
$\,w_{ij}=|(u_i,e_j)_H|^2\,$, where $\,\{e_1,e_2,\ldots\}\,$ is a
complete orthonormal subset of a separable complex Hilbert space
$\,H\,$, $\,\{u_1,u_2,\ldots\}\,$ is an orthonormal subset of the
same Hilbert space $\,H\,$ and $\,(\cdot,\cdot)_H\,$ is the inner
product in $\,H\,$.
\end{definition}

If the set $\,\{u_1,u_2,\ldots\}\,$ is also complete then the
inner products $\,(u_i,e_j)_H\,$ coincide with entries of a
unitary matrix. In this case the corresponding matrix
$\,\wb\in\UC\SC^r\,$ is doubly stochastic and is said to be {\sl
unistochastic\/}. In the finite dimensional case every
matrix $\,\wb\in\UC\SC^r\,$ is unistochastic.

\begin{definition}\label{D3.4}
If $\,\xb=\{x_1,x_2,\ldots\}\in\R^\infty\,$, let
\begin{equation}\label{3.2}
R_m^+(\xb)\
:=\sup_{\{x_{j_1},\ldots,x_{j_m}\}}\sum\limits_{n=1}^mx_{j_n}
\quad\text{and}\quad
R_m^-(\xb)\
:=\inf_{\{x_{j_1},\ldots,x_{j_m}\}}\sum_{n=1}^mx_{j_n}\,,
\end{equation}
where $\,m=1,2,\ldots\,$ and the supremum and infimum are taken
over all subsets of $\,\xb\,$ containing $\,m\,$ elements. Denote
by $\,Q_\xb\,$ the set of all sequences
$\,\yb=\{y_1,y_2,\ldots\}\in\R^\infty\,$ such that
\begin{equation}\label{3.3}
R_m^-(\xb)\ \le\ \sum_{n=1}^m y_{i_n}\ \le\ R_m^+(\xb)
\end{equation}
for each $\,m=1,2,\ldots,p$ and each collection of $\,m\,$
distinct positive integers $\,i_1,\ldots,i_m\,$. Finally, let
$\,X_{Q_\xb}\,$ be the subspace of $\,\R^\infty\,$ spanned by
$\,Q_\xb\,$.
\end{definition}

By Remark \ref{R1.11}, $\,X_{Q_\xb}\,$ is the minimal symmetric
perfect space containing $\,\xb\,$ whenever
$\,\xb\not\in\R_{00}^\infty\,$ and $\,X_{Q_\xb}=l^1\,$ whenever
$\,\xb\in l^1\setminus\{0\}\,$.

\begin{definition}\label{D3.5}
If $\,\xb=\{x_1,x_2,\ldots\}\in\R^\infty\,$, let
\begin{enumerate}
\item[]
$\,\VC_\xb\,$ be the linear space of matrices $\,\wb\,$ such that
$\,\xb\in\DC(\wb)\,$;
\item[]
$\,P_\xb^{\Gb_1}\,$, $\,S_\xb^{\Gb_1}\,$ and $\,US_\xb^r\,$ be the
sets of all sequences $\,\yb\in\R^\infty\,$ such that
$\,\yb=\wb\xb\,$ for some $\,\wb\in\PC^{\Gb_1}\,$,
$\,\wb\in\SC^{\Gb_1}\bigcap\VC_\xb\,$ and
$\,\wb\in\UC\SC^r\bigcap\VC_\xb\,$ respectively and
$S^r_\xb:=S^{\Gb_r}_\xb\,$.
\end{enumerate}
\end{definition}

Obviously, the sets $\,S^{\Gb_1}_\xb\,$, $\,P^{\Gb_1}_\xb\,$
$\,Q_\xb\,$ do not depend on the order of entries in the sequence
$\,\xb\,$. We have $\,S_\xb^\Gb\subset S_\xb^r\subset
S_\xb^\emptyset\,$, $\,P^{\Gb_r}_\xb=P^r_\xb\subset
US_\xb^r\subseteq S_\xb^r\,$ and $\,P^\Gb_\xb=P_\xb\,$ for all
$\,\xb\in\R^\infty\,$ (see Definition \ref{D1.6}).

\begin{lemma}\label{L3.6}
Let $\,\xb:=\{x_1,x_2,\ldots\}\in\R^\infty\,$,
$\,\{e_1,e_2,\ldots\}\,$ be a complete orthonormal subset of a
separable complex Hilbert space $\,H\,$ and $\,A\,$ be the
self-adjoint operator in $\,H\,$ such that $\,Ae_j=x_j\,e_j\,$.
Then $\,\yb\in US_\xb^r\,$ if and only if there exists an
orthonormal set $\,\{u_i\}\subset\DC(|A|^{1/2})\,$ such that
$\,y_i:=(Au_i,u_i)_H\,$.
\end{lemma}

\begin{proof}
A sequence $\,\yb\,$ belongs to $\,US_\xb^r\,$ if and only if
$\,y_i=\sum_j|(u_i,e_j)_H|^2x_j\,$, where $\,\{u_i\}\,$ is an
orthonormal set such that $\,\sum_j|(u_i,e_j)_H|^2|x_j|<\infty\,$
for each $\,i=1,2,\ldots\,$ These estimates are equivalent to the
inclusion $\,\{u_i\}\subset\DC(|A|^{1/2})\,$. Since
$\,u_i=\sum_j(u_i,e_j)_H\,e_j\,$, we have
$\,y_i=\sum_j|(u_i,e_j)_H|^2x_j=(Au_i,u_i)\,$.
\end{proof}

\subsection{The sets $\,P_\xb^r\,$ $\,US_\xb^r\,$, $\,S_\xb^r\,$ and $\,Q_\xb\,$}
The main result of this subsection is Theorem \ref{T3.10} which
clarifies the relation between these sets. Given a sequence
$\,\xb\,$ and a set $\,\La\subset\hat\R\,$, we shall denote by
$\,\xb\bigcap\La\,$ the sequence obtained from $\,\xb\,$ by
removing all its entries lying outside $\,\La\,$.

\begin{lemma}\label{L3.7}
Assume that the sequence $\,\xb\in\R^\infty\,$ has one
accumulation point $\,\la\in\hat\R\,$, $\,\yb\in Q_\xb\,$ and
$\,\yb\bigcap(-\infty,\la)=\emptyset\,$. Then $\,\yb\in
US_\xb^r\,$ provided that
\begin{enumerate}
\item[(a)]
either $\,\xb\bigcap[\la,+\infty)\,$ is infinite and
$\,\sharp\{i:y_i=\la\}\le\sharp\{j:x_j=\la\}\,$
\item[(b)]
or $\,\xb\bigcap[\la,+\infty)\,$ is finite and
$\,\sum_j(x_j-\la)_+ -\sum_i(y_i-\la)_+=\ve>0\,$.
\end{enumerate}
\end{lemma}

\begin{proof}
Let $\,A\,$ be defined as in Lemma \ref{L3.6}. In order to prove
the inclusion $\,\yb\in US_\xb^r\,$, we have to find an
orthonormal set $\,\{u_1,u_2,\ldots\}\subset\DC(|A|^{1/2})\,$ such
that $\,y_i:=(Au_i,u_i)_H\,$.

Assume first that (a) holds true. Then there are two entries
$\,x_{j_1},x_{k_1}\in\xb\bigcap[\la,+\infty)\,$ such that
$\,y_1\in[x_{j_1},x_{k_1}]\,$ and
$\,\xb\bigcap(x_{j_1},x_{k_1})=\emptyset\,$. If $\,y_1=\al
x_{i_1}+(1-\al)x_{k_1}\,$ and
$\,u_1:=\al^{1/2}e_{j_1}+(1-\al)^{1/2}e_{k_1}\,$ then
$\,\|u_1\|_H=1\,$ and $\,y_1=(Au_1,u_1)_H\,$. Let $\,\xb^{(1)}\,$
be the sequence obtained from $\,\xb\,$ by replacing the two
entries $\,x_{i_1}\,$ and $\, x_{k_1}\,$ with one entry
$\,x_{i_1}+ x_{k_1}-y_1\,$ and $\,\yb^{(1)}\,$ be the sequence
obtained from $\,\yb\,$ by removing the entry $\,y_1\,$. The
entries of $\,\xb^{(1)}\,$ coincide with the eigenvalues of the
self-adjoint operator $\,A_1:=\left.\Pi_1A\right|_{H_1}\,$ in the
Hilbert space $\,H_1:=\Pi_1H\,$, where $\,\Pi_1\,$ is the
orthogonal projection onto the annihilator of $\,u_1\,$.

If $\,y_1=\la\,$ then at least one of the entries
$\,x_{j_1},x_{k_1}\,$ coincides with $\,\la\,$, which implies that
$\,\xb^{(1)}\,$ and $\,\yb^{(1)}\,$ are obtained from $\,\xb\,$
and $\,\yb\,$ by removing one entry $\,\la\,$. Therefore the
sequences $\,\xb^{(1)}\,$ and $\,\yb^{(1)}\,$ satisfy the
condition (a). We also have $\,\yb^{(1)}\in Q_{\xb^{(1)}}\,$.
Indeed, if the number of entries in $\,\xb\,$ lying in the
interval $\,(x_{k_1},+\infty)\,$ is equal to $\,p\,$ then
$\,R_m^+(\xb^{(1)})=R_m^+(\xb)\,$ whenever $\,m<p\,$. If $\,m\ge
p\,$ then
$\,R_m^+(\xb^{(1)})=\,R_{m+1}^+(\xb)-y_1\ge\sum_{k=1}^my_{l_k}\,$
for each subset $\,\{y_{l_1},\ldots,y_{l_m}\}\subset\yb^{(1)}\,$.

Applying the same procedure to $\xb^{(i-1)}$, $\yb^{(i-1)}$ and
$A_{i-1}\,$ with $\,i=2,3,\ldots\,$, we can find
$\,x_{j_i},x_{k_i}\in\xb^{(i-1)}\,$, $\,u_i\in H\,$ and
$\,\xb^{(i)}\,$ such that $\,y_i\in[x_{j_i},x_{k_i}]\,$,
$\,\xb^{(i-1)}\bigcap(x_{j_i},x_{k_i})=\emptyset\,$,
$\,\Pi_{i-1}u_i=0\,$, $\,\|u_i\|_H=1\,$,
$\,y_i=(A_{i-1}u_i,u_i)_H\,$ and $\,R_m^+(\yb^{(i)})\le
R_m^+(\xb^{(i)})\,$. The entries of $\,\xb^{(i)}\,$ coincide with
the eigenvalues of $\,A_i:=\left.\Pi_iA_{i-1}\right|_{\Pi_iH}\,$,
where $\,\Pi_i\,$ is the orthogonal projection onto the
annihilator $\,H_i\,$ of the set $\,\{u_1,\ldots,u_i\}\,$. The set
$\,\{u_1,u_2,\ldots\}\,$, obtained by induction in $\,i\,$, is
orthonormal and every its element $\,u_i\,$ is a finite linear
combination of the eigenvectors $\,e_1,e_2,\ldots\,$ The latter
implies that $\,u_i\in\DC(A)\subset\DC(|A|^{1/2})\,$ and
$\,y_i=(A_{i-1}u_i,u_i)_H=(Au_i,u_i)_H\,$ for all
$\,i=1,2,3,\ldots$

If (b) holds true then $\,\la\,$ is an accumulation point of
$\,\xb\bigcap(-\infty,\la)\,$. Without loss of generality we may
assume that the sequence $\,\xb\bigcap(-\infty,\la)\,$ converges
to $\,\la\,$ and that $\,\sum_j(\la-x_j)_+<\ve/2\,$ (this can
always be achieved by removing a collection of entries from
$\,\xb\,$). Let us denote $\,\xb^{(0)}:=\xb\,$,
$\,\yb^{(0)}:=\yb\,$ and apply the same procedure as above with
$\,x_{j_i},x_{k_i}\,$ defined as follows:
\begin{enumerate}
\item[]
$\,x_{k_i}\,$ is the smallest entry of $\,\xb^{(i-1)}\,$ lying in
the interval $\,[y_i,+\infty)\,$,
\item[]
$\,x_{j_i}\,$ is either the largest entry of $\,\xb^{(i-1)}\,$
lying in $\,(\la,y_i)\,$ or, if such an entry does not exists,
$\,x_{j_i}\,$ is an arbitrary entry of
$\,\xb\bigcap(-\infty,\la)\,$.
\end{enumerate}
The inequality $\,\sum_j(\la-x_j)_+<\ve/2\,$ implies that
$\,R_m^+(\yb^{(i)})<R_m^+(\xb^{(i)})-\ve/2\,$ for all
$\,i,m=1,2,\ldots\,$. Therefore, by induction in $\,i\,$, we can
find the required representation for all entries $\,y_i\,$.
\end{proof}

\begin{lemma}\label{L3.8}
Assume that the sequence $\,\xb\in\R^\infty\,$ has two
accumulation points $\,\la,\mu\in\hat\R\,$ and $\,\la<\mu\,$. If
$\,\yb=\yb\bigcap[\la,\mu]\,$,
$\,\sharp\{i:y_i=\la\}\le\sharp\{j:x_j\le\la\}\,$ and
$\,\sharp\{i:y_i=\mu\}\le\sharp\{j:x_j\ge\mu\}\,$ then $\,\yb\in
US_\xb^r\,$.
\end{lemma}

\begin{proof}
Under the conditions of the lemma, there exists a set of distinct
positive integers $\,\{j_1,j_2,\ldots,k_1,k_2,\ldots\}\,$ such
that $\,y_i\in[x_{j_i},x_{k_i}]\,$ for all $\,i=1,2,\ldots\,$ If
$\,y_i=\al_ix_{j_i}+(1-\al_i)x_{k_i}\,$ and $\,A\,$ is defined as
in Lemma \ref{L3.6} then $\,y_i=(Au_i,u_i)_H\,$, where
$\,u_i:=\al_i^{1/2}e_{j_i}+(1-\al_i)^{1/2}e_{k_i}\,$.
\end{proof}

\begin{definition}\label{D3.9}
If $\,\xb=\{x_1,x_2,\ldots\}\in\R^\infty\,$, let
$\,x^-:=\liminf_{j\to\infty}x_j\in\hat\R\,$,
$\,x^+:=\limsup_{j\to\infty}x_j\in\hat\R\,$ and $\,\hat\xb\,$ be
the sequence obtained from $\,\xb\,$ by adding infinitely many
entries $\,x^-\,$ whenever $\,x^->-\infty\,$ and infinitely many
entries $\,x^+\,$ whenever $\,x^+<+\infty\,$.
\end{definition}

\begin{theorem}\label{T3.10} For every
$\,\xb=\{x_1,x_2,\ldots\}\in\R^\infty\,$
we have
\begin{equation}\label{3.4}
US_\xb^r\ =\ S_\xb^r\ \subseteq\ Q_\xb=Q_{\hat\xb}\ =\
US_{\hat\xb}^r\,.
\end{equation}
\end{theorem}

\begin{proof}
The equality $\,Q_\xb=Q_{\hat\xb}\,$ immediately follows from the
definition of $\,Q_\xb\,$. If $\,\yb\in S_\xb^r\,$ then for every
collection of $\,m\,$ distinct positive integers
$\,i_1,\ldots,i_m\,$ we have $\,\sum_{k=1}^my_{i_k}=\sum_j
\al_jx_j\,$, where $\,\al_j\in[0,1]\,$ and $\,\sum_j\al_j=m\,$.
This implies (\ref{3.3}). Therefore $\,S_\xb^r\subseteq Q_\xb\,$.

It remains to prove that $\,\yb\in US_\xb^r\,$ provided that
either $\,\yb\in S_\xb^r\,$ or $\,\xb=\hat\xb\,$ and $\,\yb\in
Q_\xb\,$. We are going to show that there exist countable families
of disjoint subsequences $\,\xb_n\subset\xb\,$ and
$\,\yb_n\subset\yb\,$ such that $\,\bigcup_n\xb_n=\xb\,$,
$\,\bigcup_n\yb_n=\yb\,$ and $\,\yb_n\in US_{\xb_n}^r\,$.
Obviously, this implies that $\,\yb\in US_\xb^r\,$. Given a
sequence $\,\zb\,$, in the rest of the proof we shall denote
$\,\zb^+:=\zb\bigcap(x^+,+\infty)\,$,
$\,\zb^-:=\zb\bigcap(-\infty,x^-)\,$,
$\,\zb_+:=\zb\bigcap[x^+,+\infty)\,$ and
$\,\zb_-:=\zb\bigcap(-\infty,x^-]\,$.

Assume first that $\,\xb=\hat\xb\,$. Then we can split $\,\xb\,$
into the union of three disjoint subsequences
$\,\xb_1,\xb_2,\xb_3\,$ such that $\,\xb_1=\hat\xb^+\,$,
$\,\xb_2=\hat\xb^-\,$, $\,\xb_3\,$ does not have any entries lying
outside $\,[x^-,x^+]\,$ and $\,\xb_3\,$ has infinitely many
entries $\,x^\pm\,$ whenever $\,x^\pm\,$ is finite. If
$\,\yb_1:=\yb_+\,$, $\,\yb_2:=\yb_-\,$ and
$\,\yb_0:=\yb\bigcap(x^-,x^+)\,$ then, by Lemmas \ref{L3.7} and
\ref{L3.8}, we have $\,\yb_n\in Q_{\xb_n}\,$ whenever $\,\yb\in
Q_\xb\,$. Therefore $\,Q_{\hat\xb}\subset US_{\hat\xb}^r\,$.

Assume that $\,\yb\in S_\xb^r\,$. We have to consider the
following possibilities:
\begin{enumerate}
\item[($1_+$)]
$\,\yb^+\ne\emptyset\,$, $\,\xb^+\,$ is infinite and
\begin{equation}\label{3.5}
\liminf_{m\to\infty}\left(R_m^+(\xb^+)-R_m^+(\yb^+)\right)\ =\
0\,;
\end{equation}
\item[($2_+$)]
$\,\yb^+\ne\emptyset\,$, $\,\xb^+\,$ is finite and (\ref{3.5})
holds true;
\item[($3_+$)]
$\,\yb^+\ne\emptyset\,$, $\,\xb^+\,$ is infinite and
\begin{equation}\label{3.6}
R_m^+(\xb^+)-R_m^+(\yb^+)\ge\ve>0\,,\qquad\forall m=1,2,\ldots\,;
\end{equation}
\item[($4_+$)]
$\,\xb^+\ne\emptyset\,$ is finite and (\ref{3.6}) holds true;
\item[($5_+$)]
$\,\yb^+=\emptyset\,$ and $\,\xb^+\,$ is infinite;
\item[($6'_+$)]
$\,\yb^+=\emptyset\,$ and $\,\xb^+=\emptyset\,$.
\end{enumerate}
Note that ($6'_+$) and the inclusion $\,\yb\in S_\xb^r\,$ imply
\begin{enumerate}
\item[($6_+$)]
$\,\yb^+=\emptyset\,$, $\,\xb^+=\emptyset\,$ and
$\,\sharp\{i:y_i=x^+\}\le\sharp\{j:x_j=x^+\}\,$.
\end{enumerate}
We shall say that $\,\xb\,$ and $\,\yb\,$ satisfy ($n_-$) if the
corresponding condition ($n_+$) is fulfilled for $\,-\xb\,$ and
$\,-\yb\,$.

Assume first that ($1_+$) holds true. By Lemma \ref{L3.7}, we have
$\,\yb^+\in SU_{\xb^+}^r\,$. Let
$\,\tilde\yb:=\yb\setminus\yb^+\,$ and
$\,\tilde\xb:=\xb\setminus\xb^+\,$ be the sequences obtained from
$\,\yb\,$ and $\,\xb\,$ by removing all the entries
$\,y_i\in\yb^+\,$ and $\,x_j\in\xb^+\,$ respectively. If
$\,\yb=\wb\xb\,$ and $\,\wb\in\SC^r\,$ then, in view of
(\ref{3.5}), the entry $\,w_{ij}\,$ of the matrix $\,\wb\,$ is
equal to zero whenever $\,x_j>x^+\,$ and $\,y_i\le x^+\,$.
Therefore $\,\tilde\yb=\tilde\wb\tilde\xb\,$, where
$\,\tilde\wb\in\SC^r\,$ is the matrix obtained from $\,\wb\,$ by
crossing out all the $\,i$th rows corresponding to
$\,y_i\in\yb^+\,$. If $\,\limsup_j\tilde x_j=\tilde x^+<x^+\,$ and
$\,\tilde\xb,\tilde\yb\,$ satisfy ($1_+$) then, in a similar
manner, we remove the subsequences
$\,\tilde\xb^+:=\tilde\xb\bigcap(\tilde x^+,+\infty)\,$ and
$\,\tilde\yb^+:=\tilde\yb\bigcap(\tilde x^+,+\infty)\,$. After
sufficiently (possibly, infinitely) many iterations we either
obtain two required families of disjoint subsequences $\,\xb_n\,$
and $\,\yb_n\,$ or end up with two remaining sequences satisfying
one of the conditions ($2_+$)--($6_+$). If ($1_-$) holds true then
we can apply the same procedure to the sequences $\,-\xb\,$ and
$\,-\yb\,$. Therefore it is sufficient to consider the sequences
$\,\xb\,$ and $\,\yb\,$ such that $\,\yb\in S_\xb^r\,$ and one of
the conditions ($2_\pm$)--($6_\pm$) is fulfilled.

Assume that ($2_+$) is fulfilled and $\,\yb=\wb\xb\,$, where
$\,\wb\in\SC^r\,$. If $\,\xb\,$ has finitely many entries
$\,x^+\,$, we define $\,\xb_\star:=\xb_+\,$. The condition
(\ref{3.5}) imply that the entry $\,w_{ij}\,$ of the matrix
$\,\wb\,$ is equal to zero whenever $\,x_j<x^+\,$ and $\,y_i\ge
x^+\,$. Therefore the number of entries in $\,\yb_+\,$ does not
exceed the number of entries in $\,\xb_\star\,$ and
$\,\yb_+=\wb_\star\xb_\star\,$ for some finite matrix
$\,\wb_\star\in\SC^r\,$. In the same way as in the proof of Lemma
\ref{L3.7} one can show that $\,\yb_+\in S_{\xb_\star}^r\,$. If
$\,\xb\,$ has infinitely many entries $\,x^+\,$ then we represent
$\,\xb\,$ as the union of two disjoint subsequences
$\,\tilde\xb\,$ and $\,\xb_\star\,$ such that
$\,\xb_\star=\hat\xb^+\,$,
$\,\tilde\xb\bigcap(x^+,+\infty)=\emptyset\,$ and $\,\tilde\xb\,$
contains infinitely many entries $\,x^+\,$. By Lemma \ref{L3.7},
we have $\,\yb_+\in US_{\xb_\star}^r\,$.

In both cases the sequences $\,\xb\setminus\xb_\star\,$ and
$\,\yb\setminus\yb_+\,$ satisfy ($6_+$). If ($2_-$) holds true
then, in a similar way, we can remove all the entries lying below
$\,x^-\,$. Therefore it is sufficient to prove the inclusion
$\,\yb\in US^r_\xb\,$ assuming that $\,\xb\,$ and $\,\yb\,$
satisfy (\ref{3.3}) and one of the conditions
($3_\pm$)--($6_\pm$).

If ($3_+$) is fulfilled then we choose a subsequence
$\,\xb^\star\,$ of the sequence $\,\xb^+\,$ in such a way that the
remaining sequence $\,\xb^+\setminus\xb^\star\,$ contains
infinitely many entries and $\,R_m^+(\yb^+)\le R_m^+(\xb^\star)\,$
for all $\,m=1,2,\ldots$ By Lemma \ref{L3.7}, $\,\yb_+\in
US_{\xb_\star}^r\,$ If we remove all entries $\,x_j\in\xb^\star\,$
and $\,y_i\in\yb^+\,$ then the remaining sequences
$\,\xb\setminus\xb^\star\,$ and $\,\yb\setminus\yb^+\,$ satisfy
($5_+$). Similarly, if ($3_-$) holds true then, after applying
this procedure to $\,-\xb\,$ and $\,-\yb\,$, we arrive at ($5_-$).
Therefore we can assume without loss of generality that $\,\xb\,$
and $\,\yb\,$ satisfy (\ref{3.3}) and one of the conditions
($4_\pm$)--($6_\pm$).

Let ($4_+$) be fulfilled. If $\,\yb_+=\emptyset\,$ then we simply
remove all the entries $\,x_j^+\ge x^+\,$ and arrive at ($6_+$).
Otherwise we choose a subsequence $\,\xb_\star\,$ of the sequence
$\,\xb\,$ in such a way that $\,\xb_+\subset\xb_\star\,$ and
$\,x^+\,$ is an accumulation point of both sequences
$\,\xb_\star\,$ and $\,\xb\setminus\xb_\star\,$. Lemma \ref{L3.7}
implies that $\,\yb_+\in US_{\xb_\star}^r\,$. Removing the
subsequences $\,\yb^+\,$, $\,\xb_\star\,$ and all remaining
entries $\,x_j>x^+\,$, we arrive at ($6_+$). If ($4_-$) is
fulfilled then, in a similar manner, we can remove the entries
$\,x_j\in(-\infty,x^-)\,$ and the entries
$\,y_i\in(-\infty,x^-]\,$ so that ($6_-$) holds true.

Finally, under conditions ($5_\pm$) or ($6_\pm$) the inclusion
$\,\yb\in US_\xb^r\,$ follows from Lemma \ref{L3.8}.
\end{proof}

Theorem \ref{T3.10} implies, in particular, that the set
$\,US_\xb^r\,$ is convex. Note that the set of matrices
$\,\UC\SC^r\,$ is not convex even in the finite dimensional case
(see Example \ref{E4.3}). Since the set $\,Q_\xb\,$ is
$\,\Tf_0$-closed, Theorem \ref{T3.10} also implies that
\begin{equation}\label{3.7}
\cch P_\xb^r \ \subseteq\ \overline{S_\xb^r}\ \subseteq\
Q_\xb\,,\qquad\forall\xb\in\R^\infty\,,
\end{equation}
where the closure is taken in any topology which is finer than
$\,\Tf_0\,$.

\begin{corollary}\label{C3.11}
Let $\,\Tf\,$ be an arbitrary topology on $\,X_{Q_\xb}\,$, which
is finer than $\,\Tf_0\,$ and coarser than the Mackey topology
$\,\Tf_\mrm(X_{Q_\xb},X'_{Q_\xb})\,$. Then
\begin{equation}\label{3.8}
\cch P_\xb^r\ =\ \overline{S_\xb^r}\ =\
Q_\xb\,,\qquad\forall\xb\in\R^\infty\,,
\end{equation}
where $\,\cch P_\xb^r\,$ and $\,\overline{S_\xb^r}\,$ are the
sequential $\,\Tf$-closures of the sets $\,\ch P_\xb^r\,$ and
$\,S_\xb^r\,$ respectively.
\end{corollary}

\begin{proof}
In view of (\ref{3.7}), it is sufficient to prove (\ref{3.8}) for
$\,\Tf=\Tf_\mrm(X_{Q_\xb},X'_{Q_\xb})\,$. In the rest of the prove
$\,\bar\Om\,$ denotes the sequential
$\,\Tf_\mrm(X_{Q_\xb},X'_{Q_\xb})$-closure of the set $\,\Om\in
X_{Q_\xb}\,$ and $\,\cch\Om\,$ is the sequential
$\,\Tf_\mrm(X_{Q_\xb},X'_{Q_\xb})$-closure of its convex hull.

Let $\,\VC_{\SC,\xb}:=\VC_\SC\bigcap\VC_\xb\,$, where
$\,\VC_\SC\,$ is the subspace introduced in Definition \ref{D2.8}.
By Lemma \ref{L3.1}, we have $\,\wb\xb\in X'_{Q_\xb}\,$ for all
$\,\wb\in\VC_{\SC,\xb}\,$ and, consequently,
$\,\xb'\otimes\xb\in\VC'_{\SC,\xb}\,$ for all $\,\xb'\in
X'_{Q_\xb}\,$. If $\,\xb\not\in l^\infty\,$ then
$\,X_{Q_\xb}=\R^\infty\,$ and
$\,\Tf_\mrm(X_{Q_\xb},X'_{Q_\xb})=\Tf_0\,$ is a metrizable
topology. If $\,\xb\in l^\infty\,$ then
$\,\VC_{\SC,\xb}=\VC_\SC\,$. Therefore Theorem \ref{T2.15} and
Remark \ref{R2.18} imply that $\,S_\xb^r\subseteq\cch
P_\xb^r=\overline{S_\xb^r}\,$.

Note that
\begin{enumerate}
\item[(*)] for each $\,\ve>0\,$ there exists
$\,\xb_\ve\in P_\xb^r\,$ such that $\,\hat\xb-\xb_\ve\in l^1\,$
and $\,\|\hat\xb-\xb_\ve\|_{l^1}<\ve\,$.
\end{enumerate}
Indeed, if $\,x^+<+\infty\,$ then we can always find a subsequence
$\,\{x_{j_k}\}_{k=1,2,\ldots}\,$ of $\,\xb\,$ such that the
$\,l^1$-norm of the sequence $\,\{x^+-x_{j_k}\}_{k=1,2,\ldots}\,$
is smaller than $\,\ve/6\,$. Similarly, if $\,x^->-\infty\,$ then
there exists a subsequence $\,\{x_{i_n}\}_{n=1,2,\ldots}\,$ such
that $\,i_n\ne j_k\,$ for all $\,k,n\,$ and the $\,l^1$-norm of
the sequence $\,\{x^--x_{i_n}\}_{n=1,2,\ldots}\,$ is smaller than
$\,\ve/6\,$. The required sequence $\,\xb_\ve\,$ is obtained from
$\,\hat\xb\,$ by replacing the entries $\,x^+\,$ and $\,x^-\,$
with $\,x_{j_{2k-1}}\,$ and $\,x_{i_{2n-1}}\,$ and changing the
entries $\,x_{j_k}\,$ and $\,x_{i_n}\,$ of the sequence
$\,\hat\xb\,$ to $\,x_{j_{2k}}\,$ and $\,x_{i_{2n}}\,$
respectively.

In view of Remark \ref{R1.12}, (*) implies that $\,\cch
P_{\hat\xb}^r\subseteq\cch P_\xb^r\,$. Since $\,Q_\xb\,$ is
sequentially closed, applying Theorem \ref{T3.10} and taking into
account the identity $\,\cch P_\xb^r=\overline{S_\xb^r}\,$, we
obtain $\,Q_\xb=S_{\hat\xb}^r=\cch P_{\hat\xb}^r\subseteq\cch
P_\xb^r=\overline{S_\xb^r}\subseteq Q_\xb\,$.
\end{proof}

\begin{remark}\label{R3.12}
If $\,\xb'\in X'\,$ contains a subsequence which converges to zero
and $\,\tilde\xb\in P_\xb^r\,$ then one can find $\,\xb_n\in
P_\xb\,$ such that $\,\langle\tilde\xb-\xb_n,\xb'\rangle\to0\,$ as
$\,n\to\infty\,$. This observation and the separation theorem
immediately imply that, under the conditions of Corollary
\ref{C3.11},
\begin{enumerate}
\item[(1)]
$\,\cch P_\xb=Q_\xb\,$ whenever $\,\xb\not\in l^1\,$,
\item[(2)]
$\,\cch P_\xb=Q^\star_\xb:=\{\yb\in
Q_\xb:y_1+y_2+\ldots=x_1+x_2+\ldots\}\,$ whenever $\,\xb\in l^1\,$
and $\,\Tf\,$ is the $\,l^1$-topology (indeed, if $\,\xb'\in
l^\infty\,$ separates $\,P_\xb\,$ and $\,\xb^\star\in
Q^\star_\xb\,$ and $\,c'\,$ is an accumulation point of the
sequence $\,\xb'\,$ then, by the above,
$\,\tilde\xb':=\{x'_1-c',x'_2-c',\ldots\}\,$ separates
$\,P_\xb^r\,$ and $\,\xb^\star\,$, which contradicts to Corollary
\ref{C3.11}).
\end{enumerate}
The latter result is well known (see, for example, \cite{Ma1},
Theorem 4.2), the former was proved in \cite{Ma1} for the topology
$\,\Tf\,$ generated by a symmetric norm which satisfies
(\ref{1.3}).
\end{remark}

\begin{remark}\label{R3.13}
By Corollary \ref{C3.11}, (\ref{3.8}) holds true in the Mackey
topology $\,\Tf_\mrm(l^\infty,l^1)\,$ whenever $\,\xb\in
l^\infty\,$. If, in addition, $\,x_j\to c\ne0\,$ as
$\,j\to\infty\,$ then, applying Corollary \ref{C3.11} to the
sequence $\,\tilde\xb:=\{x_1-c,x_2-c,\ldots\}\,$, one can show
that (\ref{3.8}) remains valid with respect to a stronger
topology.
\end{remark}

\subsection{Extreme points} Theorem \ref{T2.4}
suggests that $\,\ex S_\xb^{\Gb_1}\subset P_\xb^{\Gb_1}\,$. In the
next theorem we prove this inclusion only under some additional
conditions.

\begin{definition}\label{D3.14}
Denote $\,S_{\xb,(m)}^{\Gb_1}:=\{\yb^{(m)}\in\R^m:\yb\in
S_\xb^{\Gb_1}\}\,$,
$\,P_{\xb,(m)}^{\Gb_1}:=\{\yb^{(m)}\in\R^m:\yb\in
P_\xb^{\Gb_1}\}\,$, $\,S_{\xb,(m)}^r:=S_{\xb,(m)}^{\Gb_r}\,$ and
$\,P_{\xb,(m)}^r:=P_{\xb,(m)}^{\Gb_r}\,$, where $\,\yb^{(m)}\,$ is
defined as in (\ref{1.1}).
\end{definition}

Clearly, $\,S_{\xb,(\infty)}^{\Gb_1}=S_\xb^{\Gb_1}\,$ and
$\,P_{\xb,(\infty)}^{\Gb_1}=P_\xb^{\Gb_1}\,$.

\begin{theorem}\label{T3.15}
Let $\,\xb=\{x_1,x_2,\ldots\}\in\R^\infty\,$ and $\,\Gb_1\,$ be a
set of rows and columns. Assume that at least one of the following
conditions is fulfilled:
\begin{enumerate}
\item[(1)]
$\,m<\infty\,$,
\item[(2)]
$\,m=\infty\,$ and either $\,\Gb_1\subseteq\Gb_r\,$ or $\,\Gb_1\,$
contains all columns,
\item[(3)]
$\,m=\infty\,$ and $\,x_j\ne x_k\,$ whenever $\,j\ne k\,$.
\end{enumerate}
Then $\,\ex S_{\xb,(m)}^{\Gb_1}\subset P_{\xb,(m)}^{\Gb_1}\,$. If
{\rm(3)} holds true and $\,\xb\,$ does not contain zero entries
then $\,\wb\xb\not\in\ex S_\xb^{\Gb_1}\,$ whenever
$\,\wb\in\left(\SC^{\Gb_1}\setminus\PC^{\Gb_1}\right)\bigcap\VC_\xb\,$.
\end{theorem}

\begin{proof}
Let $\,\wb\in\SC^{\Gb_1}\bigcap\VC_\xb\,$ and
$\,(\wb\xb)^{(m)}\in\ex S_{\xb,(m)}^{\Gb_1}\,$. The proof consists
of two parts. In the first part we shall construct a special
matrix $\,\tilde\wb\in\SC^{\Gb_1}\bigcap\VC_\xb\,$ such that
$\,(\tilde\wb\xb)^{(m)}=(\wb\xb)^{(m)}\,$. Then we shall show that
$\,(\tilde\wb\xb)^{(m)}=(\wb_0\xb)^{(m)}\,$ with some
$\,\wb_0\in\PC^{\Gb_1}\,$ and that
$\,\wb=\tilde\wb\in\PC^{\Gb_1}\,$ whenever (3) holds true and
$\,\xb\,$ does not have zero entries.

Let $\,\La\,$ be the countable set of all distinct values taken by
the entries of $\xb\,$, $\,\Jb_\la=\{j_1,j_2,\ldots\}\,$ be the
ordered set of all indices $\,j_1<j_2<\dots\,$ such that
$\,x_{j_k}=\la\,$ and
$\,v_{i,\la}:=\sum_{j_k\in\Jb_\la}w_{ij_k}\,$.

If (3) is fulfilled then we take $\,\tilde\wb:=\wb\,$. Otherwise,
given $\,m\le+\infty\,$ and an ordered set
$\,\Jb_\la=\{j_1,j_2,\ldots\}\,$, we define
\begin{enumerate}
\item[(1)]
$\,\tilde w^{(m;\la)}_{1j_1}:=v_{1,\la}\,$ and $\,\tilde
w^{(m;\la)}_{1j_k}:=0\,$ for all $k>1$;
\item[(2)]
if $\,1<i\le m\,$ and $\,j_l\in\Jb_\la\,$ is the maximal positive
integer such that $\,\tilde w_{(i-1)j_l}>0\,$, then
\begin{enumerate}
\item[]
$\,\tilde w^{(m;\la)}_{ij_k}:=0\,$ for all $\,k<l\,$ and
$\,k>l+1\,$,
\item[]
$\,\tilde w^{(m;\la)}_{ij_l}:=\min\{\,v_{i,\la}\,,
1-\sum_{n=1}^{i-1}\tilde w^{(m;\la)}_{nj_l}\,\}\,$ and
\item[]
$\,\tilde w^{(m;\la)}_{ij_{l+1}}:=v_{i,\la}-w^{(m;\la)}_{ij_l}\,$.
\end{enumerate}
\end{enumerate}
Let $\,\tilde\wb^{(m)}\,$ be the $\,m\times\infty$-matrix whose
entries $\,\tilde w^{(m)}_{ij}\,$ coincide with  $\,\tilde
w^{(m;\la)}_{ij_k}\,$ for all $\,j=j_k\in\Jb_\la\,$. Obviously, we
have $\,\tilde\wb^{(m)}\in\VC_\xb\,$ and
$\,\tilde\wb^{(m)}\xb=(\wb\xb)^{(m)}\,$. For each $\,\la\in\La\,$
the matrix $\,\tilde\wb^{(m)}\,$ has at most two nonzero entries
lying at the intersections of a given $i$th row and the
$\,\Jb_\la$-columns. The minimal column-number of such a nonzero
entry in the $(i+1)$th row is not smaller than the maximal
column-number of a nonzero entry in the $i$th row; in other words,
the set of nonzero entries lying in the $\,\Jb_\la$-columns is
ladder-shaped.

If (2) is fulfilled then we take $\,\tilde\wb:=\tilde\wb^{(m)}\,$.
The matrix $\,\tilde\wb^{(m)}\,$ has the same row-sums as
$\,\wb\,$ and its column-sums $\,u_j:=\sum_{i=1}^m\tilde
w^{(m)}_{ij}\,$ are not greater than 1. If all column-sums of
$\,\wb\,$ are equal to 1 then each column-sum of
$\,\tilde\wb^{(m)}\,$ is also equal to 1. Therefore
$\,\tilde\wb\in\SC^{\Gb_1}\,$.

If (1) is fulfilled, let us consider the ordered set
$\,\Jb=\{j_1,j_2,\ldots\}\,$ of all indices $\,j_1<j_2<\dots\,$
such that $\,u_{j_k}<1\,$. Note that every set $\,\Jb_\la\,$
contains at most one element of $\,\Jb\,$. Let
$\,\tilde\wb=\{\tilde w_{ij}\}\,$ be the
$\,\infty\times\infty$-extension of the $\,m\times\infty$-matrix
$\,\tilde\wb^{(m)}\,$ defined as follows:
\begin{enumerate}
\item[(1)]
$\,\tilde w_{ij}=0\,$ for all $\,j\not\in\Jb\,$
and $\,i>m\,$;
\item[(2)]
$\,\tilde w_{(m+1)j_1}:=1-u_{j_1}\,$ and $\,\tilde
w_{(m+1)j_k}:=\min\{1-u_{j_k}\,,\,1-\sum_{n=1}^{k-1}\tilde
w_{(m+1)j_n}\}\,$ for all $\,j_k\in\Jb\,$ with $\,k=2,3,\ldots$;
\item[(3)]
if $\,i>m+1\,$ and $\,j_l\in\Jb\,$ is the maximal positive integer
such that $\,\tilde w_{(i-1)j_l}>0\,$ then
\begin{enumerate}
\item[]
$\,\tilde w_{ij_k}:=0\,$ for all $\,j_k\in\Jb\,$ with $k<l$,
\item[]
$\,\tilde w_{ij_l}:=1-u_{j_l}-\tilde w_{(i-1)j_l}\,$ and
\item[]
$\,\tilde w_{ij_k}:=\min\{1-u_{j_k}\,,\,1-\sum_{n=1}^{k-1}\tilde
w_{ij_n}\}\,$ for all $\,j_k\in\Jb\,$ with $\,k>l\,$.
\end{enumerate}
\end{enumerate}
We have $\,\sum_{j\in\Jb}\,u_j\le m\,$ and, consequently,
$\,\sum_{j\in\Jb}\,(1-u_j)=+\infty\,$. Therefore, for each
$\,i>m\,$, the set of nonzero entries in the $\,i$th row of the
matrix $\,\tilde\wb\,$ is finite and non-empty. Since
$\,\tilde\wb^{(m)}\in\VC_\xb\,$, this implies that
$\,\tilde\wb\in\VC_\xb\,$. All column-sums of the matrix
$\,\tilde\wb\,$ are equal to 1. Its $\,i$th row-sum coincides with
the $\,i$th row sum of $\,\wb\,$ whenever $\,i\le m\,$ and is equal to 1
whenever $\,i>m\,$. Therefore
$\,\tilde\wb\in\SC^{\Gb_1}\,$. The set of nonzero entries
$\,\tilde w_{ij}\,$ with $\,i>m\,$ is also ladder-shaped. More
precisely, the $j$th column contains at most two such nonzero
entries (if it does then these entries lie in adjacent rows)
and the minimal column-number of a nonzero entry in the $(i+1)$th
row is not smaller than the maximal column-number of a nonzero
entry in the $i$th row.

Let $\,\tilde G\,$ the subgraph of $G$, which contains all the
vertices $g_{ij}$ (that is, the intersections of $\,i$th rows and
$\,j$th columns) such that $\,\tilde w(g_{ij}):=\tilde
w_{ij}\in(0,1)\,$. Denote $\,\tilde G_\la:=\{g_{ij}\in\tilde
G:j\in\Jb_\la\}\,$ and $\,\tilde G':=\{g_{ij}\in\tilde
G:i>m\}\,$.

Assume first that $\,\tilde G\,$ contains an admissible cycle
$\,g_{i_1j_1}\to g_{i_2j_1}\to g_{i_2j_2}\to\dots\to
g_{i_1j_1}\,$. Replacing $\,g_{i_kj_k}\to g_{i_{k+1}j_k}\to
g_{i_{k+1}j_{k+1}}\to\dots\to g_{i_{k+l}j_{k+l}}\,$ with
$\,g_{i_kj_k}\to g_{i_{k+l}j_{k+l}}\,$ whenever $\,i_k=i_{k+l}\,$,
we obtain an admissible cycle $\,g_{i'_1j'_1}\to g_{i'_2j'_1}\to
g_{i'_2j'_2}\to\dots\to g_{i'_1j'_1}\,$ which has at most two
vertices in every row. By our construction, the subgraphs
$\,\tilde G_\la\,$ and $\,\tilde G'\,$ are ladder-shaped and, for
every $\,\la\in\La\,$, the intersection $\,\tilde
G_\la\bigcap\tilde G'\,$ contains at most one element. Therefore
this admissible cycle has at least two vertices lying in the same
$\,i$th row with $\,i\le m\,$ but in distinct sets $\,\tilde
G_\la\,$. If $\,\wb_\ve^\pm\,$ are defined as in the part ({\bf2})
of the proof of Theorem \ref{T2.4} then
$\,\wb_\ve^\pm\in\SC^{\Gb_1}\bigcap\VC_\xb\,$,
$\,\tilde\wb=\frac12\,(\wb_\ve^++\wb_\ve^-)\,$ and
$\,(\wb_\ve^+\xb)^{(m)}\ne(\wb_\ve^-\xb)^{(m)}\,$. Therefore
$\,(\tilde\wb\xb)^{(m)}\not\in\ex S_{\xb,(m)}^{\Gb_1}\,$.

Thus, the graph $\,\tilde G\,$ does not have any admissible
cycles. Let us take an arbitrary vertex $\,g_{ij_0}=g_0\in\tilde
G\,$ with $\,i\le m\,$, define $\,\wb_\ve^\pm\,$ as in the part
({\bf3}) of the proof of Theorem \ref{T2.4} and denote
$\,\wb^*:=\frac12\,(\wb_\ve^+-\wb_\ve^-)\,$. Then
$\,\wb_\ve^\pm\in\SC^{\Gb_1}\bigcap\VC_\xb\,$,
$\,\tilde\wb=\frac12\,(\wb_\ve^++\wb_\ve^-)\,$ and
$$
w_{ij_0}^\star=\pm\ve\,\tilde w_{ij_0}\,,\qquad
w_{ij}^\star=\mp\ve\,\tilde w_{ij_0}\,(1-\tilde w_{ij_0})^{-1}\,\tilde
w_{ij}\,,\quad\forall j\ne j_0\,.
$$
Since $\,(\tilde\wb\xb)^{(m)}\in\ex S_{\xb,(m)}^{\Gb_1}\,$, we
have $\,\wb^\star\,\xb=0\,$ which implies that $\,\tilde
w_{ij_0}x_{j_0}=\,\tilde w_{ij_0}\,(1-\tilde
w_{ij_0})^{-1}\sum_{j\ne j_0}\tilde w_{ij}x_j\,$ and
$\,x_{j_0}=\sum_{j=1}^\infty \tilde w_{ij}x_j\,$. The integer
$\,j_0\,$ can be chosen in an arbitrary way. Therefore for each
$\,i\le m\,$ we have either $\,\tilde w_{ij}=0\,$ or
$\,j\in\Jb_{\la_i}\,$, where $\,\la_i:=\sum_{j=1}^\infty \tilde
w_{ij}x_j\in\La\,$. The first row of $\,\tilde\wb\,$ may contain
only one nonzero entry $\,w_{1j}\,$ with $\,j\in\Jb_{\la_1}\,$. If
it does then $\,x_j=w_{1j}x_j\,$ and, consequently, either
$\,w_{1j}=1\,$ or $\,x_j=\la_1=0\,$. By induction in $\,i\,$, the
same is true for all $\,i=1,2,\ldots,m\,$; namely, each of the
first $\,m\,$ rows either contains one entry 1 in a
$\,\Jb_{\la_i}$-column corresponding to some $\,\la_i\ne0\,$ or
has nonzero entries only in the $\,\Jb_0$-columns. If
$\,\Jb_0=\emptyset\,$, this implies that
$\,\tilde\wb\in\PC^{\Gb_1}\,$. If $\,\Jb_0\ne\emptyset\,$ then, by
Remark \ref{R2.16}, there exists a matrix
$\,\wb_0\in\PC^{\Gb_1}\,$ whose entries at the intersections of
the first $\,m\,$ rows and the $\,\Jb_{\la_i}$-columns with
$\,\la_i\ne0\,$ coincide with the corresponding entries of
$\,\tilde\wb\,$. Since
$\,(\wb\xb)^{(m)}=(\tilde\wb\xb)^{(m)}=(\wb_0\xb)^{(m)}$, this
completes the proof.
\end{proof}

\section{Applications to spectral theory}

\subsection{Notation and definitions}\label{S4.1}
Let $H$ be a separable complex Hilbert space with the inner
product $(\cdot,\cdot)_H$ and norm $\|\cdot\|_H\,$. For the sake
of definiteness, we shall be assuming that $\dim H=\infty$; the
finite dimensional versions of our results are either well known
or can be proved in a similar manner.

Consider a linear operator $A$ in $H$ and denote by
$\,Q_A[\cdot]\,$ its quadratic form defined on the domain
$\DC(Q_A):=\DC(|A|^{1/2})\,$. We shall always be assuming that the
operator $\,A\,$ is self-adjoint. Let $\,\si(A)\,$,
$\,\si_\crm(A)\,$, and $\,\si_\ess(A)\,$ be its spectrum,
continuous spectrum and essential spectrum respectively and let
$\,\si_\prm(A)=\{\la_1,\la_2,\ldots\}\,$ be the set of its
eigenvalues. As usual, we enumerate the eigenvalues $\,\la_j\,$
taking into account their multiplicities. If $\,\La\in\hat\R\,$
and $\,\R\bigcap\La\,$ is a Borel set, we shall denote by
$\,\Pi_\La\,$ and $\,A_\La\,$ the spectral projection of $\,A\,$
corresponding to $\,\R\bigcap\La\,$ and the restriction of $\,A\,$
to the subspace $\,\Pi_\La H\,$ respectively.

\begin{definition}\label{D4.1}
Let $\,\hat\si^\pm_\ess(A)\,$ and $\,\hat\si_\ess(A)\,$ be the
subsets of $\,\hat\R\,$ such that
\begin{enumerate}
\item[]
$\la\in\hat\si^+_\ess(A)\,$ if and only if
$\,\dim\Pi_{[\la,\mu)}H=\infty\,$ for all $\,\mu>\la\,$,
\item[]
$\la\in\hat\si^-_\ess(A)\,$ if and only if
$\,\dim\Pi_{(\mu,\la]}H=\infty\,$ for all $\,\mu<\la\,$,
\end{enumerate}
and $\hat\si_\ess(A):=\hat\si^-_\ess(A)\bigcup\hat\si^+_\ess(A)\,$
\end{definition}

Obviously, $\,\si_\ess(A)=\R\bigcap\hat\si_\ess(A)\,$,
$\,+\infty\not\in\si^+_\ess(A)\,$ and
$\,-\infty\not\in\si^-_\ess(A)\,$. We have
$\,\pm\infty\in\hat\si_\ess(A)\,$ if and only if $\,\pm A\,$ is
not bounded from above.

\begin{definition}\label{D4.2}
If $\,m\,$ is a positive integer or $\,m=\infty\,$, let
\begin{enumerate}
\item[(1)]
$\si(m,A)\,$ be the set of vectors
$\,\xb=(x_1,x_2,\ldots)\in\R^m\,$ such that $\,x_j\in\si(A)\,$ for
each $\,j\,$ and the number of entries
$\,x_j=\la\not\in\si_\ess(A)\,$ does not exceed the multiplicity
of the eigenvalue $\,\la\,$;
\item[(2)]
$\si_\prm(m,A)\,$ be the set of vectors
$\,\xb=(x_1,x_2,\ldots)\in\R^m\,$ such that
$\,x_j\in\si_\prm(A)\,$ for each $\,j\,$ and the number of entries
$\,x_j=\la\,$ does not exceed the multiplicity of the eigenvalue
$\,\la\,$.
\end{enumerate}
If $\,\ub=\{u_1,u_2,\ldots\}\,$ is an orthonormal subset of
$\DC(Q_A)$ which contains $m$ elements $\,u_k\,$, denote
$\,Q_A[\ub]:=\{Q_A[u_1],Q_A[u_2],\ldots\}\in\R^m\,$ and define
\begin{enumerate}
\item[(3)]
$\Si(m,A)\ :=\ \{\,\yb\in\R^m:\yb=Q_A[\ub]\,$ for some
$\,\ub\subset\DC(Q_A)\,\}\,$.
\end{enumerate}
The sets $\,\si(m,A)\,$, $\,\si_\prm(m,A)\,$ and $\,\Si(m,A)\,$
will be called the {\sl $m$-spectrum}, {\sl point $m$-spectrum},
and {\sl $m$-numerical range} of $A$ respectively.
\end{definition}

The $m$-spectra and $m$-numerical range are symmetric with respect
to permutations of the coordinates $x_k$. The $\infty$-spectra and
$\infty$-numerical range are subsets of $\R^\infty$, whose
projections onto the subspace spanned by any $m$ coordinate
vectors coincide with the $\,m$-spectra and $\,m$-numerical range.
In particular, $\,\si(1,A)=\si(A)\,$,
$\,\si_\prm(1,A)=\si_\prm(A)\,$ and $\,\Si(1,A)\,$ is the
numerical range of the operator $\,A\,$. Since $\si(A)$ is a
closed set, the $m$-spectrum $\,\si(m,A)\,$ is closed in the
topology of element-wise convergence and, consequently, in any
finer topology.

Definition \ref{D4.2} can be extended to an arbitrary linear
operator $\,A\,$ acting in the separable Hilbert space $\,H\,$. In
\cite{BD}, Section 36, the authors defined a matrix $m$-numerical
range as the set of all $m\times m$-matrices of the form $\Pi
A\Pi$, where $\Pi$ is an orthogonal projection of rank $m<\infty$.
Halmos defined an $m$-numerical range as the set of traces of such
matrices (see \cite{H}, Chapter 17). Our definition lies in
between: we consider the sets of diagonal elements of the matrices
$\Pi A\Pi$ instead of their traces. Yet another concept of
multidimensional numerical range, related to a given block
representation of the operator $\,A\,$, was introduced in
\cite{LMMT}. Halmos' $m$-numerical range is always convex (in the
self-adjoint case this immediately follows from Corollary
\ref{C4.7}). The $m$-numerical range $\Si(m,A)$ is convex if $A$
is self-adjoint. The matrix $\,m$-numerical range considered in
\cite{BD} and the multidimensional numerical range introduced in
\cite{LMMT} are not necessarily convex. The latter depends on the
choice of block representation and is not unitary invariant.

If $\,A\ne A^*\,$ then  $\Si(m,A)\,$ does not have to be convex,
even if the operator $A$ is normal and $\dim H<\infty\,$. The
following simple example was suggested by A. Markus \cite{Ma2}.

\begin{example}\label{E4.3}
Let $\,A=\{a_{ij}\}\,$ be the diagonal $3\times3$-matrix with
$\,a_{11}=i\,$, $\,a_{22}=1\,$ and $\,a_{33}=0\,$. Then
$\,\{i,1,0\}\in\Si(3,A)\,$ and $\,\{0,i,1\}\in\Si(3,A)\,$.
However, the half-sum $\,\{\frac i2,\frac 12+\frac i2,\frac12\}\,$
does not belong to $\,\Si(3,A)\,$. In the same way as in Lemma
\ref{L3.6}, one can show that $\,\Si(3,A)=\bigcup_\wb\wb\,\zb\,$,
where $\,\zb\,$ is the three dimensional complex vector
$\,\{0,i,1\}\,$ and the union is taken over all unistochastic
$\,3\times3$-matrices $\,\wb\,$. This implies that the set of
unistochastic matrices is not convex.
\end{example}

In \cite{FW} the authors proved that
$\cch\si_\prm(m,A)=\cch\Si(m,A)$ whenever $A$ is a normal $m\times
m$-matrix. There are also some results on the so-called
$c$-numerical range of a finite matrix $A$, which is defined as
the image of $\Si(m,A)$ under the map
$\xb\to\langle\xb,c\rangle\in\C$ where $c$ is a fixed
$m$-dimensional complex vector (see \cite{GR}, \cite{MMF},
\cite{MS}).

\subsection{Extreme points of the multidimensional numerical
range} We shall need the following simple lemma.

\begin{lemma}\label{L4.4}
If $\,\si(A)\subset[\la_-,\la_+]\,$ and
$\,\la_\pm\in\hat\si_\ess(A)\,$ then $\,\Si(\infty,A)\,$ coincides
with the set of all sequences $\,\zb=\{z_1,z_2,\ldots\}\,$ such
that $\,z_i\in[\la_-,\la_+]\,$ for all $\,i\,$ and the number of
entries $\,z_i=\la_\pm\,$ does not exceed the multiplicity of the
eigenvalue $\,\la_\pm\,$ (we assume that the multiplicity is zero
whenever $\,\la_\pm\,$ is not an eigenvalue).
\end{lemma}

\begin{proof}
The spectral theorem implies that every sequence
$\,\zb\in\Si(\infty,A)\,$ satisfies the above two conditions. On
the other hand, if $\,z_1\in[\la_-,\la_+]\,$ then, using the
spectral theorem, one can easily find $\,u_1\in\DC(A)\,$ such that
$\,\|u\|_H=1\,$ and $\,z_1=Q_A[u_1]\,$. Clearly, $\,u_1\,$ is an
eigenvector whenever $\,z_1=\la_-\,$ or $\,z_1=\la_+\,$. If
$\,\Pi_1\,$ is the orthogonal projection onto the annihilator of
$\,u_1\,$  and $\,A_1:=\Pi_1A\Pi_1\,$ then $\,\DC(A_1)=\DC(A)\,$
and $\,A-A_1\,$ is a finite rank operator. Since a finite rank
perturbation does not change the essential spectrum, by induction
in $\,i\,$ we can construct an orthonormal set
$\,\ub=\{u_i\}\subset\DC(A)\,$ such that $\,z_i=Q_A[u_i]\,$ for
all $\,i=1,2,\ldots$
\end{proof}

\begin{definition}\label{D4.5}
We shall say that $\,\xb=\{x_1,x_2,\ldots\}\in\si(\infty,A)\,$ is
a {\sl generating\/} sequence of the self-adjoint operator $\,A\,$
if
\begin{enumerate}
\item[(1)]
either $\,\si_\crm(A)=\emptyset\,$,
$\,\xb\subset\si_\prm(\infty,A)\,$ and $\,\xb\,$ contains all the
eigenvalues $\,\la_j\,$ of $\,A\,$ according to their
multiplicities;
\item[(2)]
or $\,\si_\crm(A)\ne\emptyset\,$ and $\,\xb\,$ can be represented
as the union of three disjoint subsequences, one of which is
defined as above and the other two lie in the open interval
$\,(\inf\si_\crm(A),\sup\si_\crm(A))\,$ and converge to
$\,\inf\si_\crm(A)\,$ and $\,\sup\si_\crm(A)\,$ respectively.
\end{enumerate}
\end{definition}

\begin{theorem}\label{T4.6}
If $\,\xb\,$ is a generating sequence of $\,A\,$ then
$\,\Si(m,A)=S_{\xb,(m)}^r\,$.
\end{theorem}

\begin{proof}
Since $\,\Si(m,A)\,$ and $\,S_{\xb,(m)}^r\,$ coincide with the
projections of $\,\Si(\infty,A)\,$ and $\,S_\xb^r\,$ onto the
subspace spanned by the first $\,m\,$ coordinate vectors, it is
sufficient to prove that $\,\Si(\infty,A)=S_\xb^r\,$. If
$\,\tilde\xb\,$ is another generating sequence then, by Lemma
\ref{L3.8}, $\,\tilde\xb\in S_\xb^r\,$. Therefore $\,S_\xb^r\,$
does not depend on the choice of generating sequence $\,\xb\,$.

Let $\,\la^-:=\inf\si_\crm(A)\,$, $\,\la^+:=\sup\si_\crm(A)\,$,
$\,\La:=(\inf\si_\crm(A),\sup\si_\crm(A))\,$, $\,\la_j\,$ be the
eigenvalues of $\,A\,$ lying outside $\,\La\,$ and $\,\{e_j\}\,$
be the orthonormal set of eigenvectors corresponding to
$\,\la_j\,$.

Assume first that $\,\yb=Q_A[\ub]\,$, where
$\,\ub\subset\DC(Q_A)\,$ is an orthonormal set. Let
$\,d_i:=\|\Pi_\La u_i\|_H\,$ and $\,\{z_i\}\subset\La\,$ be a
sequence with two accumulation points $\,\la^\pm\,$, such that
$\,Q_A[\Pi_\La u_i]=d_i^2(\al_i\,z_{2i-1}+(1-\al_i)\,z_{2i})\,$
with some $\,\al_i\in[0,1]\,$. Then
$$
Q_A[u_i]=Q_A[\Pi_\La u_i]+Q_A[\Pi_{\R\setminus\La}u_i]=
d_i^2\,\al_i\,z_{2i-1}+d_i^2\,(1-\al_i)\,z_{2i}+\sum_jw_{ij}\la_j\,,
$$
where $\,w_{ij}:=|(u_i,e_j)_H|^2\,$. Since
$\,\sum_jw_{ij}=\|\Pi_{\R\setminus\La}u_i\|_H^2=1-d_i^2\,$ and
$\,\sum_iw_{ij}\le\|e_j\|_H^2\le1\,$, this implies that $\,\yb\in
S_\xb^r\,$, where $\,\xb\,$ is an arbitrary generating sequence
containing all the eigenvalues $\,\la_j\,$ and the subsequence
$\,\{z_i\}\,$.

Assume now that $\,\yb=\{y_1,y_2,\ldots\}\in S_\xb^r\,$ for some
generating sequence $\,\xb\,$. By Lemma \ref{L4.4}, there exists
an orthonormal set $\,\{v_n\}\subset \Pi_\La H\bigcap\DC(Q_A)\,$
such that $\,x_n=Q_A[v_n]\,$ for all $x_n\in\La\,$. Let $\,\tilde
A\,$ be the self-adjoint operator in the space $\,H\,$ such that
$\,\tilde Ae_j=\la_je_j\,$ and $\,\tilde Av_n=x_n\,v_n\,$ for all
$x_n\in\La\,$. In view of Lemma \ref{L3.6} and Theorem
\ref{T3.10}, we have $\,Q_{\tilde A}[\tilde\ub]=\yb\,$ for some
orthonormal set $\,\tilde\ub=\{\tilde u_1,\tilde
u_2,\ldots\}\subset\DC(Q_{\tilde A})\,$. If $\,\tilde
d_i:=\|\Pi_\La \tilde u_i\|_H\,$ and $\,\tilde z_i:=\tilde
d_i^{-2}Q_{\tilde A}[\Pi_\La\tilde u_i]\,$ then the sequence
$\,\{\tilde z_i\}\,$ satisfies conditions of Lemma \ref{L4.4}.
Therefore $\,\tilde z_i=Q_A[u'_i]\,$ for some orthonormal set
$\{u'_i\}\subset\Pi_\La H\bigcap\DC(Q_A)\,$. Since
$\,A_{\R\setminus\La}=\tilde A_{\R\setminus\La}\,$, the
orthonormal set $\,\ub:=\{d_iu'_i+\Pi_{\R\setminus\La}\tilde
u_i\}\,$ satisfies $\,Q_A[\ub]=\yb\,$.
\end{proof}

\begin{corollary}\label{C4.7}
For each $\,m=1,2,\ldots,\infty\,$ the set $\,\Si(m,A)\,$ is
convex and $\,\ex\Si(m,A)\subset\si_\prm(m,A)\,$. A sequence
$\,\yb\in\si_\prm(m,A)\,$ belongs to $\,\ex\Si(m,A)\,$ if and only
if there is a (possibly, degenerate) interval
$\,[\mu^-,\mu^+]\subset\hat\R\,$ such that
\begin{enumerate}
\item[(1)]
$\,\si_\crm(A)\subset[\mu^-,\mu^+]\,$,
$\,\hat\si^+_\ess(A)\bigcap[-\infty,\mu^-)=\emptyset\,$,
$\,\hat\si^-_\ess(A)\bigcap(\mu^+,+\infty]=\emptyset\,$;
\item[(2)]
$\,\yb\bigcap(\mu^-,\mu^+)=\emptyset\,$ and $\,\yb\,$ contains all
the eigenvalues $\,\la_j\not\in[\mu^-,\mu^+]\,$ according to their
multiplicities.
\end{enumerate}
\end{corollary}

\begin{proof}
Let $\,\xb\,$ be a generating sequence. Theorems \ref{T3.15} and
\ref{T4.6} imply that the set $\,\Si(m,A)=S_{\xb,(m)}^r\,$ is
convex and $\,\ex\Si(m,A)=\ex S_{\xb,(m)}^r\subset
P_{\xb,(m)}^r\,$.

Let $\,\yb\in\si_\prm(m,A)\,$ and $\,\mu^\pm\in\hat\R\,$ satisfy
(1) and (2). Then $\,\yb\bigcap(-\infty,\mu^-]\,$ either is empty
or coincides with the union of disjoint nondecreasing subsequences
$\,\yb_n\,$ such that $\,\sup\yb_n\le\inf\yb_{n+1}\,$ and
$\,\sup\yb_n\not\in\yb_n\,$ whenever $\,\yb_n\,$ is infinite (in
the latter case $\,A\,$ is bounded from below). Using this
observation, one can easily show by induction in $\,n\,$ that the
sequence $\,\yb\bigcap(-\infty,\mu^-]\,$ cannot be represented as
a convex combination of two distinct sequences from
$\,S_{\xb,(k)}^r\,$. Similarly, $\,\yb\bigcap[\mu^+,+\infty)\,$ is
not a convex combination of two distinct sequences from
$\,S_{\xb,(k)}^r\,$. Therefore every sequence
$\,\yb\in\si_\prm(m,A)\,$ satisfying the conditions of the
corollary belongs to $\,\ex S_{\xb,(m)}^r\,$.

Assume now that $\,\yb\in\ex S_{\xb,(m)}^r\,$ and denote
$\,\si_\yb:=\{\la\in\si(A):\la\not\in\yb\}\,$. If
$\,\hat\si_\ess^-(A)\ne\emptyset\,$,
$\,\hat\si_\ess^+(A)\ne\emptyset\,$,
$\,\inf\hat\si_\ess^+(A)<\sup\hat\si_\ess^-(A)\,$ and
$$
\yb_\star\ :=\
\yb\bigcap(\inf\hat\si_\ess^+(A),\sup\hat\si_\ess^-(A))\ \ne\
\emptyset
$$
then $\,\yb_\star\,$ coincides with a convex combination of two
distinct sequences $\,\yb_\star^\pm\,$ whose entries lie in the
open interval $\,(\inf\hat\si_\ess^+(A),\sup\hat\si_\ess^-(A))\,$.
By Lemma \ref{L3.8}, we have $\,\yb_\star^\pm\in
S_{\xb_\star,(k)}^r\,$, where $\,k\,$ is the number of entries in
$\,\yb_0\,$ and
$\,\xb_\star:=\xb\bigcap(\inf\hat\si_\ess^+(A),\sup\hat\si_\ess^-(A))\,$.
Therefore the sequence $\,\yb\in\ex S_{\xb,(m)}^r\,$ does not have
entries which are greater than $\,\inf\hat\si_\ess^+(A)\,$ and
smaller than $\,\sup\hat\si_\ess^-(A)\,$. In particular,
$\,\yb\bigcap(\inf\si_\crm(A),\sup\si_\crm(A))=\emptyset\,$. Since
the number of entries $\,\inf\si_\crm(A)\,$ and
$\,\sup\si_\crm(A)\,$ in the generating sequence $\,\xb\,$ does
not exceed the multiplicity of the corresponding eigenvalue and
$\,\yb\in P_{\xb,(m)}^r\,$, this implies that
$\,\yb\in\si_\prm(m,A)\,$.

Let
\begin{enumerate}
\item[]
$\,\mu^-=\mu^+:=\inf\hat\si_\ess^+(A)\,$ if
$\,\si_\yb=\emptyset\,$ and $\,\hat\si_\ess^-(A)=\emptyset\,$;
\item[]
$\,\mu^-=\mu^+:=\sup\hat\si_\ess^-(A)\,$ if
$\,\si_\yb=\emptyset\,$ and $\,\hat\si_\ess^+(A)=\emptyset\,$;
\item[]
$\,\mu^-=\mu^+:=\mu\,$ if $\,\si_\yb=\emptyset\,$ and
$\,\inf\hat\si_\ess^+(A)\ge\sup\hat\si_\ess^-(A)\,$, where
$\,\mu\,$ is an arbitrary number from the closed interval
$\,[\sup\hat\si_\ess^-(A),\inf\hat\si_\ess^+(A)]\,$;
\item[]
$\,\mu^-:=\inf\si_\yb\,$ and $\,\mu^+:=\sup\si_\yb\,$ if
$\,\si_\yb\ne\emptyset\,$.
\end{enumerate}
Obviously, in the first three cases (1) and (2) hold true. It
remains to prove that $\,\yb\bigcap(\mu^-,\mu^+)=\emptyset\,$,
$\,\hat\si^+_\ess(A)\bigcap(-\infty,\mu^-)=\emptyset\,$ and
$\,\hat\si^-_\ess(A)\bigcap(\mu^+,+\infty)=\emptyset\,$ in the
last case.

Let $\,\si_\yb\ne\emptyset\,$ and $\,\mu^\pm\,$ be defined as
above. If $\,\si_\yb\,$ contains two distinct entries $\,\la\,$
and $\, \mu\,$ and $\,\yb\,$ has an entry $\,y_i\in(\la,\mu)\,$
then $\,\yb\,$ coincides with a convex combination of two distinct
sequences obtained by replacing $\,y_i\,$ with $\,\la\,$ and
$\,\mu\,$ respectively. Both these sequences belong to
$\,P_{\xb,(m)}^r\,$ for some generating sequence $\,\xb\,$.
Therefore the inclusion $\,\yb\in\ex S_{\xb,(m)}^r\,$ implies that
$\,\yb\bigcap(\mu^-,\mu^+)=\emptyset\,$.

If $\,m<\infty\,$ then
$\,\hat\si^+_\ess(A)\bigcap(-\infty,\mu^-)=\emptyset\,$ as the
number of eigenvalues lying below $\,\mu^-\,$ is finite. Assume
that $\,m=\infty\,$ and that there exists
$\,\hat\la\in\hat\si^+_\ess(A)\,$ such that $\,\la<\mu^-\,$. Let
$\,\yb^\star\,$ be a decreasing subsequence of
$\,\yb\bigcap(-\infty,\mu^-)\,$, which converges to $\,\la\,$, and
$\,\xb^\star\in P_\xb^r\,$ be the sequence obtained from
$\,\yb^\star\,$ by adding an entry $\,\mu\in\si_\yb\,$. By Lemma
\ref{L3.7}, the sequences $\,\yb_\pm^\star\,$ obtained from
$\,\yb^\star\,$ by replacing an arbitrary entry
$\,y_i\in\yb^\star\,$ with $\,y_i-\ve>\la\,$ and
$\,y_i+\ve<\mu^-\,$ respectively belong to $\,S_{\xb^*}^r\,$.
Therefore $\,\yb_\pm\in S_\xb^r\,$, where $\,\yb_\pm\,$ are the
sequences obtained from $\,\yb\,$ by replacing the entry $\,y_i\,$
with $\,y_i\pm\ve\,$. Since $\,\yb=\frac12(\yb_-+\yb_+)\,$, this
contradicts to the inclusion $\,\yb\in\ex S_\xb^r\,$.

In a similar way one can show that
$\,\hat\si^-_\ess(A)\bigcap(\mu^+,+\infty)=\emptyset\,$.
\end{proof}

\begin{remark}\label{R4.8}
Let $\,\La_\erm(A)\subset\hat\R\,$ be the intersection of all
intervals $\,[\mu^-,\mu^+]\,$ satisfying the condition (1) of
Corollary \ref{C4.7}. If the number of eigenvalues lying outside
$\,\La_\erm(A)\,$ is smaller then $\,m\,$ then, by Corollary
\ref{C4.7}, the set $\,\Si(m,A)\,$ does not have any extreme
points.
\end{remark}

\begin{definition}\label{D4.9}
If $\,\xb\,$ is a generating sequence of $\,A\,$, let
$\,Q(\infty,A):=Q_\xb\,$, $\,X_A:=X_{Q_\xb}\,$  (see Definition
\ref{D3.5}), $\,Q(m,A)\,$ be the projection of $\,Q(\infty,A)\,$
on the subspace of $\,X_A\,$ spanned by the first $\,m\,$
coordinate vectors and $\,\Tf_A^{(m)}\,$ be the topology on
$\,Q_A^{(m)}\,$ induced by $\,\Tf_\mrm(X_A,X'_A)\,$.
\end{definition}

Obviously, the symmetric perfect space $\,X_A\,$ and its subset
$\,Q(\infty,A)\,$ do not depend on the choice of generating
sequence $\,\xb\,$. By Theorems \ref{T3.10} and \ref{T4.6}, we
have
\begin{equation}\label{4.1}
S_\xb^r\ \subseteq\ \Si(\infty,A)\ \subseteq\ Q(\infty,A)
\end{equation}
for each generating sequence $\,\xb\,$.

\begin{lemma}\label{L4.10}
For every $\,\xb\in\si(\infty,A)\,$ and every sequence of strictly
positive numbers $\ve_k$ there exists $\,\yb\in\Si(\infty,A)\,$
such that $\,|y_k-x_k|\le\ve_k\,$. For every
$\,\yb\in\Si(\infty,A)\,$ there exists a sequence of vectors
$\,\yb_n\in\ch\si(\infty,A)\,$ which converges to $\,\yb\,$ in the
Mackey topology $\,\Tf_\mrm(X_A,X'_A)\,$.
\end{lemma}

\begin{proof}
Let $\,\La_k:=(x_k-\ve_k,x_k+\ve_k)\,$. If $\,x_j\in\si_\ess(A)\,$
then $\,\dim P_{\La_j}H=\infty\,$. Since a finite dimensional
perturbation does not change the essential spectrum, by induction
in $\,k\,$ one can find an orthonormal sequence
$\,\{u_1,u_2,\ldots\}\,$ such that $\,Au_k=x_ku_k\,$ whenever
$\,x_k\not\in\si_\ess(A)\,$ and $\,u_k\in P_{\La_j}H\,$ otherwise.
If $\,y_k=Q_A[u_k]\,$ then $\,\yb\in\Si(\infty,A)\,$ and
$\,|y_k-x_k|\le\ve_k\,$.

The second statement of the lemma follows from Corollary
\ref{C3.11} and the second inclusion (\ref{4.1}).
\end{proof}

Lemma \ref{L4.10} immediately implies that
\begin{equation}\label{4.2}
\cch\si(m,A)\ = \overline{\Si(m,A)}\ =\ Q(m,A)\,, \qquad\forall
m=1,2,\ldots,\infty\,,
\end{equation}
where the bar denotes the sequential closure taken in any
topology which is finer than the topology of element-wise
convergence $\,\Tf_0$ and coarser than $\,\Tf_A^{(m)}\,$. Since
$\,\Tf_0\,$ is a metrizable topology, (\ref{4.2}) remains valid if
we take the usual closure.

\begin{corollary}\label{C4.11}
For each $\,m=1,2,\ldots,\infty\,$ the set $\,Q(m,A)\,$ is convex
and $\,\ex Q(m,A)\subset\si(m,A)\,$. A sequence
$\,\yb\in\si(m,A)\,$ belongs to $\,\ex Q(m,A)\,$ if and only if
there is a (possibly, degenerate) interval
$\,[\mu^-,\mu^+]\subset\hat\R\,$ such that
\begin{enumerate}
\item[(1)]
$\,\hat\si_\ess(A)\subset[\mu^-,\mu^+]\,$,
\item[(2)]
$\,\yb\bigcap(\mu^-,\mu^+)=\emptyset\,$ and $\,\yb\,$ contains all
the eigenvalues $\,\la_j\not\in[\mu^-,\mu^+]\,$ according to their
multiplicities.
\end{enumerate}
\end{corollary}

\begin{proof}
Let $\,\xb\,$ be a generating sequence, $\,x^+=\limsup\xb\,$ and
$\,x^-=\liminf\xb\,$. In view of (\ref{3.4}), (\ref{3.8}) and
(\ref{4.2}), we have $\,Q(m,A)=S_{\hat\xb,(m)}^r\,$. Therefore the
corollary is obtained by applying Corollary \ref{C4.7} to the
operator $\,A\oplus A_+\oplus A_-\,$ acting in the orthogonal sum
$\,H\oplus H_+\oplus H_-\,$, where  $\,A_\pm\,$ is multiplication
by $\,x^\pm\,$ in $\,H_\pm\,$, $\,\dim H_\pm=\infty\,$ whenever
$\,|x^\pm|<\infty\,$ and $\,H_\pm=\emptyset\,$ otherwise.
\end{proof}

\begin{remark}\label{R4.12}
By Corollary \ref{C4.11}, each sequence $\,\yb\in\ex Q(m,A)\,$
consists of eigenvalues $\,\la_j\not\in\cch\si_\ess(A)\,$  and,
possibly, a collection of entries $\,\inf\si_\ess(A)\,$ and
$\,\sup\si_\ess(A)\,$. All these eigenvalues can be found with the
use of the Rayleigh--Ritz variational formula. The interval
$\,\La_\erm(A)\,$ defined in Remark \ref{R4.8} is a subset of
$\,\hat\si_\ess(A)\,$ and may be strictly smaller. Therefore a
sequence $\,\yb\in\ex\Si(\infty,A)\,$ may contain eigenvalues
lying inside $\,\cch\si_\ess(A)\,$.
\end{remark}

\begin{example}\label{E4.13}
Assume that the continuous spectrum of $\,A\,$ is empty and that
the eigenvalues of $\,A\,$ form a sequence $\,\xb\,$ which has two
accumulation points $\,\la^\pm\,$ such that $\,\la^+>\la_-\,$.
Then $\,\cch\si_\ess(A)=[\la^-,\la^+]\,$. However, if $\,\la^-\,$
or $\,\la^+\,$ is not an accumulation point of the sequence
$\,\xb\bigcap[\la^-,\la^+]\,$ then $\,\La_\erm(A)=\emptyset\,$ and
$\,\xb\,$ is an extreme point of $\,\Si(\infty,A)\,$.
\end{example}

\begin{example}\label{E4.14}
If $\,\hat\si_\ess(A)=[-\infty,+\infty]\,$ then
$\,\Si(\infty,A)=Q(\infty,A)=\R^\infty\,$ and
$\,\ex\Si(\infty,A)=\emptyset\,$. If
$\,\hat\si_\ess(A)=\{+\infty\}\,$ then
$\,\Si(\infty,A)=Q(\infty,A)\,$ and the extreme points of
$\,\Si(\infty,A)\,$ are the sequences formed by all the
eigenvalues $\,\la_j\,$. If $\,\hat\si_\ess(A)=[\mu,+\infty]\,$,
$\,\mu\in\hat\si^+_\ess(A)\bigcap\R\,$ and $\,\xb\,$ is the
sequence formed by all the eigenvalues $\,\la_j<\mu\,$ then every
extreme point of $\,Q(\infty,A)\,$ is obtained from $\,\xb\,$ by
adding an arbitrary collection of entries $\,\mu\,$ and every
extreme point of $\,\Si(\infty,A)\,$ is obtained from $\,\xb\,$ by
adding a collection of entries $\,\mu\,$ whose number does not
exceed the multiplicity of the eigenvalue $\,\mu\,$ (we assume
that the multiplicity is zero if $\,\mu\,$ is not an eigenvalue).
\end{example}

\begin{remark}\label{R4.15}
Let $\,\R^\infty_\la\,$ be the set of all real sequences with
entries in the interval $\,(-\infty,\la]\,$. Theorem \ref{T4.6}
implies that
$\,\Si(\infty,A)\bigcap\R^\infty_\la=\Si(\infty,A_{(-\infty,\la]})\,$
whenever $\,\rank A_{(\la,+\infty)}=\infty\,$. This observation
allows one to extend Theorem \ref{T4.6} and Corollaries
\ref{C4.7}, \ref{C4.11} to the sets
$\,\Si(\infty,A)\bigcap\R^\infty_\la\,$ and
$\,\si(\infty,A)\bigcap\R^\infty_\la\,$. Note that the linear
space $\,X_{A_\la-\la\Irm}\,$ may well be smaller than $\,X_A\,$.
In this case one can refine Lemma \ref{L4.10} and related results
by considering the operator $\,A_\la-\la\Irm\,$ instead of
$\,A\,$.
\end{remark}

\subsection{Variational formulae and exposed points}
Recall that a function $\,\ps:\Om\to\hat\R\,$ defined on a convex
set $\,\Om\,$ is called {\sl quasi-concave\/} if
\begin{equation}\label{4.3}
\ps(\al\,\xb+(1-\al)\,\yb)\ \ge\ \min\{\ps(\xb),\ps(\yb)\}\,,
\qquad\forall\xb,\yb\in\Om\,,\ \forall\al\in(0,1)\,,
\end{equation}
and {\sl strictly quasi-concave\/} if the left hand side of
(\ref{4.3}) is strictly greater than the right hand side. The
function $\ps$ is quasi-concave if and only if the sets
$\,\{\xb\in X:\ps(\xb)\ge\la\}\,$ are convex for all
$\la\in\hat\R$. The function $\,\ps\,$ is said to be sequentially
upper $\,\Tf$-semicontinuous if these sets are sequentially closed
in the topology $\,\Tf\,$. The identity (\ref{4.2}) and Corollary
\ref{C4.7} immediately imply the following two variational
results.

\begin{corollary}\label{C4.16}
If $\,\ps\,$ is a quasi-concave sequentially upper
$\,\Tf_A^{(m)}$-semi\-conti\-nuous function on $\,Q(m,A)\,$ then
\begin{equation}\label{4.4}
\inf_{\xb\in\si(m,A)}\ps(\xb)\ =\ \inf_{\xb\in\Si(m,A)}\ps(\xb)\,.
\end{equation}
\end{corollary}

For each finite $m$ the functions $\,\ps(\xb)= x_1+x_2\dots+x_m\,$
and $\,\ps(\xb)=x_1\,x_2\ldots x_m=\exp(\ln x_1+\dots\ln x_m)\,$
defined on the set of positive sequences are quasi-concave and
$\,\Tf_0$-upper semicontinuous. Therefore the variational formulae
for the sum and product of the first $m$ eigenvalues of a positive
self-adjoint operator are particular cases of (\ref{4.4}).

\begin{corollary}\label{C4.17}
Let $\,\ps\,$ be a real-valued function defined on $\,\Si(m,A)\,$.
If
\begin{enumerate}
\item[(a)]
either $\,\ps\,$ is quasi-concave and
$\,\ps(\yb)<\ps(\tilde\yb)\,$ for all $\,\tilde\yb\ne\yb\,$
\item[(b)]
or $\,\ps\,$ is strictly quasi-concave and
$\,\ps(\yb)\le\ps(\tilde\yb)\,$ for all $\,\tilde\yb\,$
\end{enumerate}
then $\,\yb\in\si_\prm(m,A)\,$.
\end{corollary}

Note that $\,\yb\,$ is a $\,\Tf$-exposed point of the set
$\,\Si(m,A)\,$ if and only if there exists a linear
$\,\Tf$-continuous function $\,\ps\,$ satisfying the condition
(a).

\begin{example}\label{E4.18}
If $\,m<\infty\,$ then $\,Q(m,A)\,$ is a closed convex polytope,
$\,\Si(m,A)\,$ is a convex dense subset of $\,Q(m,A)\,$ and, by
Corollaries \ref{C4.7} and \ref{C4.11}, we have
$\,\ex\Si(m,A)\subset\ex Q(m,A)\,$. In this case the extreme
points of $\,\Si(m,A)\,$ and $\,Q(m,A)\,$ are exposed in the
standard Euclidean topology.
\end{example}

The sets $\,\Si(\infty,A)\,$ and $\,Q(\infty,A)\,$ may contain
extreme points which are not $\,\Tf_\mrm(X_A,X'_A)$-exposed.

\begin{example}\label{E4.19}
If $\,A\,$ is not bounded then $\,X_A=\R^\infty\,$ and
$\,X'_A=\R_{00}^\infty\,$. For every $\,\yb\in\R^\infty\,$ and
$\,\xb'\in\R_{00}^\infty\,$ there exists $\,\tilde\yb\in P_\yb\,$
such that $\,\tilde\yb\ne\yb\,$ and
$\,\langle\yb,\xb'\rangle=\langle\tilde\yb,\xb'\rangle\,$.
Therefore the sets $\,\Si(\infty,A)\,$ and $\,Q(\infty,A)\,$ do
not contain $\,\Tf_\mrm(X_A,X'_A)$-exposed points whenever $\,A\,$
is unbounded.
\end{example}

If $\,\yb\in\ex\Si(\infty,A)\,$ or $\,\yb\in\ex Q(\infty,A)\,$,
let $\,[\mu^-,\mu^+]\,$ be the interval introduced in Corollary
\ref{C4.7} or \ref{C4.11} respectively,
$\,\yb_{(+)}:=\yb\bigcap[\mu^+,+\infty)\,$,
$\,\yb_{(-)}:=\yb\bigcap(-\infty,\mu^-]\,$ and $\,\la^\pm\,$ be
defined as follows:
\begin{enumerate}
\item[]
$\,\la_\yb^+:=\limsup\yb_{(+)}\,$ whenever $\,\yb_{(+)}\,$ is
infinite, $\,\la_\yb^+:=\inf\yb_{(+)}\,$ whenever $\,\yb_{(+)}\,$
is finite and nonempty, and $\,\la_\yb^+:=\mu^+\,$ whenever
$\,\yb_{(+)}=\emptyset\,$;
\item[]
$\,\la_\yb^-:=\liminf\yb_{(-)}\,$ whenever $\,\yb_{(-)}\,$ is
infinite, $\,\la_\yb^-:=\sup\yb_{(-)}\,$ whenever $\,\yb_{(-)}\,$
is finite and nonempty, and $\,\la_\yb^-:=\mu^-\,$ whenever
$\,\yb_{(-)}=\emptyset\,$.
\end{enumerate}
If $\,\la_\yb^-<\la_\yb^+\,$, denote by $\,\La_\yb\,$ the interval
with end points $\,\la_\yb^-\,$ and $\,\la_\yb^+\,$ such that
$\,\la_\yb^\pm\in\La_\yb\,$ if and only if $\,\la_\yb^\pm\,$ is an
accumulation point of the sequence obtained from $\,\yb\,$ by
removing all the entries $\,y_j\in[\la_\yb^-,\la_\yb^+]\,$. If
$\,\la_\yb^-=\la_\yb^+\,$, let
$\,\La_\yb:=[\la_\yb^-,\la_\yb^-]\,$.

Obviously, $\,\si_\ess(A)\subset\bar\La_\yb\,$ and $\,\yb\,$
contains all the eigenvalues lying outside the closure
$\,\bar\La_\yb\,$ of the interval $\,\La_\yb\,$. The entries of
$\,\yb\,$ lying below and above $\,\La_\yb\,$ form a nondecreasing
sequence $\,\yb^{(-)}\,$ and a nonincreasing sequence
$\,\yb^{(+)}\,$ respectively (either of these sequences may be
empty).

\begin{theorem}\label{T4.20}
If $\,A\,$ belongs to the trace class then every extreme point
$\,\yb\in\ex Q(\infty,A)\,$ or $\,\yb\in\ex\Si(\infty,A)\,$ is
$\,\Tf_\mrm(X_A,X'_A)$-exposed. If $\,A\,$ is bounded but does not
belong to the trace class then
\begin{enumerate}
\item[]
$\,\yb\in\ex Q(\infty,A)\,$ is a $\,\Tf_\mrm(X_A,X'_A)$-exposed
point of $\,Q(\infty,A)\,$ if and only if either
$\,\yb\bigcap\La_\yb=\emptyset\,$ or $\,\La_\yb\,$ consists of one
point;
\item[]
$\,\yb\in\ex\Si(\infty,A)\,$ is a $\,\Tf_\mrm(X_A,X'_A)$-exposed
point of $\,\Si(\infty,A)\,$ if and only if either
$\,\yb\bigcap\La_\yb=\emptyset\,$ or $\,\La_\yb\,$ is closed and
the spectrum of the truncation $\,A_{\La_\yb}\,$ consists of one
point.
\end{enumerate}
\end{theorem}

\begin{proof}
Assume that $\,\tilde\yb\in Q(\infty,A)\,$ or
$\,\yb\in\ex\Si(\infty,A)\,$ and
$\,\tilde\yb\in\Si(\infty,A)\subset Q(\infty,A)\,$. Let
$\,y_{j_1}\le y_{j_2}\le\ldots\,$ be the entries of
$\,\yb^{(-)}\,$, $\,y_{k_1}\ge y_{k_2}\ge\ldots\,$ be the entries
of $\,\yb^{(+)}\,$ and $\,y_{n_1},y_{n_2},\ldots\,$ be the entries
of $\,\yb\,$ lying in $\,\La_\yb\,$. Consider an arbitrary
sequence $\,\xb'\in X'_A\,$ such that
$$
x'_{j_1}<x'_{j_2}<\ldots<0\,,\quad
x'_{k_1}>x'_{k_2}>\ldots>0\quad\text{and}\quad
x'_{n_1}=x'_{n_2}=\ldots=0\,.
$$
The identity (\ref{1.6}) implies that
$\,\sum_iy_{j_i}\,x'_{j_i}\ge\sum_i\tilde y_{j_i}\,x'_{j_i}\,$ and
these two sums coincide only if $\,\tilde y_{j_i}=y_{j_i}\,$ for
all $\,i\,$. Similarly, $\,\sum_iy_{k_i}\,x'_{j_i}\ge\sum_i\tilde
y_{k_i}\,x'_{k_i}\,$ and the sums coincide only if $\,\tilde
y_{k_i}=y_{k_i}\,$ for all $\,i\,$. If  $\,\La_\yb\,$ satisfies
the conditions of the theorem, $\,\tilde y_{j_i}=y_{j_i}\,$ for
all $\,i\,$ and $\,\tilde y_{k_i}=y_{k_i}\,$ for all $\,i\,$ then,
in view of Theorem \ref{T4.6}, we have $\,\tilde\yb=\yb\,$.
Therefore the sequence $\,\yb\,$ is
$\,\Tf_\mrm(X_A,X'_A)$-exposed.

If $\,\yb\in l^1\,$ and
$$
x'_{j_1}<x'_{j_2}<\ldots<-2\,,\quad
x'_{k_1}>x'_{k_2}>\ldots>2\,,\quad x'_{n_1}=x'_{n_2}=\ldots=1
$$
then the same arguments show that
$\,\langle\yb,\xb'\rangle>\langle\tilde\yb,\xb'\rangle\,$ for all
$\,\tilde\yb\in Q(\infty,A)\,$. This proves the first statement of
the theorem.

Assume now that $\,A\,$ does not belong to the trace class and
that $\,\yb\,$ is $\,\Tf_\mrm(X_A,X'_A)$-exposed. Then there
exists a sequence $\,\xb'\in X'_A\subset\R_0^\infty\,$ such that
$\,\langle\yb,\xb'\rangle>\langle\tilde\yb,\xb'\rangle\,$ whenever
$\,\tilde\yb\in P_\yb\,$ and $\,\tilde\yb\ne\yb\,$. If
$\,y_i>y_j\,$ but $\,x'_i\le x'_j\,$ then
$\,\langle\yb,\xb'\rangle\le\langle\tilde\yb,\xb'\rangle\,$, where
$\,\tilde\yb\in P_\yb\,$ is the sequence obtained from $\,\yb\,$
by interchanging the entries $\,y_i\,$ and $\,y_j\,$. Therefore
\begin{enumerate}
\item[(c$_2$)]
$\,x'_i>x'_j\,$ whenever $\,y_i>y_j\,$.
\end{enumerate}

If $\,\la_\yb^-=\la_\yb^+\,$ then $\,\La_\yb\,$ satisfies the
conditions of the theorem. Assume that $\,\la_\yb^-<\la_\yb^+\,$.
Then $\,\la_\yb^\pm\,$ are accumulation points of $\,\yb\,$ and
$\,\La\,$ is not closed if and only if $\,\yb\,$ contains
infinitely many entries $\,\la_\yb^-\,$ or $\,\la_\yb^+\,$. The
inclusion $\,\xb'\in\R_0^\infty\,$ and (c$_2$) imply that
$\,x_i=0\,$ whenever $\,y_i\in\La_\yb\,$. If $\,\yb\,$ has two
distinct entries in $\,\La_\yb\,$ then
$\,\langle\yb,\xb'\rangle=\langle\tilde\yb,\xb'\rangle\,$, where
$\,\tilde\yb\ne\yb\,$ is the sequence obtained by interchanging
these entries. Therefore either $\,\yb\bigcap\La_\yb=\emptyset\,$
or there exists $\,\la\,$ such that $\,y_i=\la\,$ whenever
$\,y_i\in\La_\yb\,$. If $\,\yb\bigcap\La_\yb\ne\emptyset\,$ and
$\,\si(A_{\La_\yb})\,$ contains another point $\,\mu\ne\la\,$ then
we can find $\,u\in\Pi_{[\la,\mu]}H\,$ such that
$\,\tilde\la:=Q_A[u]\ne\la\,$ and the sequence $\,\tilde\yb\,$
obtained by replacing $\,\la\,$ with $\,\tilde\la\,$ belongs to
$\,\Si(\infty,A)\,$. Since
$\,\langle\yb,\xb'\rangle=\langle\tilde\yb,\xb'\rangle\,$, we see
that $\,\si(A_{\La_\yb})=\{\la\}\,$ whenever
$\,\yb\bigcap\La_\yb\ne\emptyset\,$. Finally, if $\,\yb\in
Q(\infty,A)\,$ or $\,\La_\yb\,$ is not closed then
$\,\langle\yb,\xb'\rangle=\langle\tilde\yb,\xb'\rangle\,$ for the
sequence $\,\tilde\yb\,$ obtained by replacing $\,\la\,$ with
$\,\la_\yb^-\,$ or $\,\la_\yb^+\,$. Therefore in either case
$\,\yb\bigcap\La_\yb=\emptyset\,$.
\end{proof}

\subsection{Family of operators}
Finally, let us consider a family of self-adjoint operators
$\,\{A_\th\}_{\th\in\The}\,$ acting in $\,H\,$, where $\,\The\,$
is an arbitrary index set. The following
corollary implies that
\begin{equation}\label{4.6}
\si(\infty,A)\ \subset\ \cch\bigcup_{\th\in\The}\si(\infty,A_\th)
\end{equation}
whenever $\,A\in\cch\{A_\th\}\,$, provided that the closures are
taken in appropriate topologies.

\begin{corollary}\label{C4.22}
Let $\,X\,$ be a subspace of $\,\R^\infty\,$ and $\,A\,$ be a
self-adjoint operator in $\,H\,$ such that $\,X_A\subset X\,$ and
$\,X_{A_\th}\subset X\,$ for all $\,\th\in\The\,$. Assume that for
every orthonormal set $\,\ub\subset\DC(Q_A)\,$, every $\,\xb'\in
X'\,$ and every $\,\ve>0\,$ there exist an operator $\,A_\th\,$
and an orthonormal set $\,\tilde\ub\subset\DC(Q_{A_\th})\,$ such
that $\,\langle Q_A[\ub],\xb'\rangle\le\langle
Q_{A_\th}[\tilde\ub],\xb'\rangle+\ve\,$. Then we have
{\rm(\ref{4.6})}, where the closure is taken in the Mackey
topology $\,\Tf_\mrm(X,X')\,$.
\end{corollary}

\begin{proof}
By the separation theorem, under conditions of the corollary we
have
$\,\Si(\infty,A)\subset\cch\bigcup_{\th\in\The}\Si(\infty,A_\th)\,$.
Therefore (\ref{4.6}) follows from (\ref{4.2}).
\end{proof}

In Corollary \ref{C4.22} we can always take $\,X=\R^\infty\,$, in
which case $\,X'=\R_{00}^\infty\,$ and $\,\Tf_\mrm(X,X')\,$
coincides with the topology of element-wise convergence
$\,\Tf_0\,$. If $\,A\,$ and $\,A_\th\,$ satisfy the conditions of
Corollary \ref{C4.22} and are compact then we can take
$\,X=\R_0^\infty\,$, which implies (\ref{4.6}) with the closure
taken in the $\,l^\infty$-topology.

\end{document}